\documentclass[10pt]{amsart}
\usepackage{amsmath,amssymb,amsthm,amscd,mathrsfs}
\usepackage{latexsym,amsmath,amssymb,graphicx}
\usepackage[all]{xy}
\usepackage{color}
\usepackage{epsf}
\xyoption{v2}
\xyoption{2cell}
\xyoption{dvips}
\CompileMatrices
\numberwithin{equation}{section}
\newtheorem{prop}{Proposition}[section]

\newtheorem{rema}[prop]{Remark}

\newtheorem{conj}[prop]{Conjecture}

\numberwithin{equation}{section}

\newcommand{\be}{\begin{equation}}
\newcommand{\ee}{\end{equation}}

\newcommand{\IP}{\mathbb{P}}

\newcommand\IZ{\mathbb {Z}}

\newcommand\IQ{\mathbb {Q}}
\newcommand{\IC}{\mathbb{C}}

\newcommand{\IR}{\mathbb{R}}

\newcommand{\ba}{\begin{array}}
\newcommand{\ea}{\end{array}}

\newcommand{\CX}{{\mathcal X}}

\newcommand{\wH}{{\widetilde H}}

\newcommand{\CY}{{\mathcal Y}}
\newcommand{\uY}{{\underline Y}}

\newcommand{\bal}{\begin{aligned}}
\newcommand{\eal}{\end{aligned}}

\newcommand{\CZ}{{\mathcal Z}}

\newcommand{\mfm}{{\mathfrak{M}}}

\newcommand{\longto}{\longrightarrow}

\newcommand{\CO}{{\mathcal O}}

\newcommand{\CE}{{\mathcal E}}

\newcommand{\CC}{{\mathcal C}}

\newcommand{\obj}{{\mathfrak {Ob}}}

\newcommand{\wgamma}{{\widetilde \gamma}}
\newcommand{\wA}{{\widetilde A}}
\CompileMatrices
\LaTeXdiagrams
\UseAllTwocells
\newdimen\tableauside\tableauside=1.0ex
\newdimen\tableaurule\tableaurule=0.4pt
\newdimen\tableaustep
\def\phantomhrule#1{\hbox{\vbox to0pt{\hrule height\tableaurule width#1\vss}}}
\def\phantomvrule#1{\vbox{\hbox to0pt{\vrule width\tableaurule height#1\hss}}}
\def\sqr{\vbox{%
  \phantomhrule\tableaustep
  \hbox{\phantomvrule\tableaustep\kern\tableaustep\phantomvrule\tableaustep}%
  \hbox{\vbox{\phantomhrule\tableauside}\kern-\tableaurule}}}
\def\squares#1{\hbox{\count0=#1\noindent\loop\sqr
  \advance\count0 by-1 \ifnum\count0>0\repeat}}
\def\tableau#1{\vcenter{\offinterlineskip
  \tableaustep=\tableauside\advance\tableaustep by-\tableaurule
  \kern\normallineskip\hbox
    {\kern\normallineskip\vbox
      {\gettableau#1 0 }%
     \kern\normallineskip\kern\tableaurule}%
  \kern\normallineskip\kern\tableaurule}}
\def\gettableau#1 {\ifnum#1=0\let\next=\null\else
  \squares{#1}\let\next=\gettableau\fi\next}
\title{Wallcrossing and Cohomology of The Moduli Space of Hitchin Pairs}
\author{Wu-yen Chuang, Duiliu-Emanuel Diaconescu, Guang Pan}

\begin{document}

\begin{abstract}
A conjectural recursive relation for the Poincar\'e polynomial of the  Hitchin moduli space is derived from wallcrossing in the refined local Donaldson-Thomas theory of a
a curve. A doubly refined generalization of this theory is also
conjectured and shown to similarly determine the Hodge polynomial of the same
moduli space.
\end{abstract}

\maketitle

\tableofcontents

\section{Introduction}
Let $X$ be a smooth projective curve over $\IC$ and $M_1,M_2$ be line bundles
on $X$ so that $M_1\otimes_X M_2\simeq K_X^{-1}$.  Any such triple
$\CX=(X,M_1,M_2)$ determines an abelian category $\CC_\CX$ of twisted quiver
sheaves on $X$, called ADHM sheaves. This construction is briefly explained in
section (\ref{review}).  ADHM sheaves are essentially $(M_1,M_2)$-twisted
representations of an ADHM quiver in ${\mathfrak{Coh}}(X)$, the $\CO_X$-module
associated to the framing node being isomorphic to $\CO_X^{\oplus v}$, for some
$v\in \IZ_{\geq 0}$. In particular $v=0$ objects of $\CC_\CX$ are Higgs sheaves
on $X$ i.e. coherent $\CO_X$-modules $E$ decorated by a morphism
$E\otimes_X(M_1\oplus M_2)\to E$ satisfying the standard integrability condition
\eqref{eq:Higgsrelation}.
Objects with $v\geq 1$ include in addition framing
data consisting of morphisms $E\otimes_X M_1\otimes_X M_2\to \CO_X^{\oplus v}$,
$\CO_X^{\oplus v} \to E$ satisfying a modified integrability condition
\eqref{eq:ADHMrelation}.

The purpose of this paper is to present an application of ADHM sheaves to
computations of Betti and Hodge numbers of moduli spaces of stable Hitchin
pairs on the curve $X$. As a brief history of the subject, note that the Poincar\'e
polynomial of the moduli space of stable bundles on a curve
  has been recursively
computed in \cite{Poincarepol_stabbundles}, \cite{cohgroups} using number theoretic
methods, respectively \cite{YM-surface} using gauge theoretic methods. The Hodge
polynomial of the same moduli spaces, has been
recursively computed in \cite{earl-kirwan}, and also in \cite{motive-bundles,hodge-pairs,
hodge-twotwo,hodge-rankthree} for bundles of rank two and three\footnote{According to \cite{earl-kirwan}, the Hodge polynomial of the moduli space of rank three bundles has
been first derived by P. Newstead in unpublished work.}.
The Poincar\'e polynomial of the moduli space of stable Hitchin pairs with coprime
rank and degree has been computed by Hithchin in \cite{hitchin-selfd} for rank two,
and Gothen,
\cite{Bettinumbers}, for rank three. Using number theoretic techniques, a  conjectural formula for any rank has been
derived by Hausel and Rodriguez-Villegas in  \cite{HRV} and generalized to Hodge polynomials by Hausel in \cite{Mirror-Hodge}. Similar results for parabolic rank three
Higgs bundles have been obtained in \cite{betti-parabolic}.
Finally, the motive of the moduli space of
rank four Hitchin pairs in the Grothedieck ring of algebraic varieties is computed in the
upcomig work \cite{motive-higgs}.

The present paper presents a string theoretic perspective on this subject based on wallcrossing and  refined generalized Donaldson-Thomas invariants.
There are currently two theories of Donaldson-Thomas invariants, the Kontsevich-Soibelman
theory \cite{wallcrossing} and the Joyce-Song
theory \cite{genDTI}. The former is based on a construction of motivic Donaldson-Thomas
invariants which specialize to integral valued invariants in a semiclassical limit.
 The later
constructs $\IQ$-valued generalized Donaldson-Thomas invariants which are
conjecturally related to these integral invariants by a multicover formula
\cite[Sect. 6.2]{genDTI}. The application presented below relies on the motivic Donaldson-Thomas
theory of Kontsevich and Soibelman applied to ADHM sheaves, or, equivalently, on a conjectural
refinement of Joyce-Song theory.

The generalized Donaldson-Thomas theory of ADHM sheaves has been studied
using the fomalism of \cite{genDTI} in
\cite{chamberI,chamberII,higherrank}.
Moduli spaces of ADHM sheaves have been constructed in \cite{chamberI}
using a natural stability
condition depending on a real parameter $\delta \in \IR$.
The main results for $v=1$ objects, which is the relevant case in this paper,
are reviewed in section (\ref{review}).
In particular  for fixed numerical invariants  $\gamma=(r,e)\in \IZ_{\geq 1}\times \IZ$ there is a finite set of critical stability parameters dividing
the real axis into stability chambers. Note that $\delta=0$ is
critical for all $(r,e)\in  \IZ_{\geq 1}\times \IZ$.
Residual ADHM invariants $A_\delta(r,e)$
are defined in each chamber
by equivariant virtual integration \cite{chamberI}. The asymptotic invariants
$A_{+\infty}(r,e)$ corresponding to $\delta>>0$ are identified with
local stable pair invariants in \cite{modADHM}.
Wallcrossing formulas for ADHM invariants
are derived in \cite{chamberII} using the formalism of Joyce \cite{J-I,J-II,J-III,J-IV} and
Joyce and Song \cite{genDTI}. The resulting wallcrossing formulas are also shown to be
in agreement with the Kontsevich-Soibelman formula \cite{wallcrossing}.
Note that the theory of Joyce and Song also yields residual generalized Donaldson-Thomas
invariants $H(r,e)$ counting semistable Higgs sheaves on $X$ with numerical invariants
$(r,e)\in \IZ_{\geq 1}\times \IZ$. These invariants enter the wallcrossing formulas
for $A_\delta(r,e)$ derived in \cite{chamberII}.

The conjectures formulated in section (\ref{motivicconjectures}) below
summarize the main results of Kontsevich-Soibelman theory needed in this paper.
Since the virtual enumerative theory of ADHM sheaves has been studied in
\cite{chamberI,chamberII} employing Joyce-Song theory, these conjectures can be
also viewed as a refinement of their generalized Donaldson-Thomas
formalism.
In particular, the conjectural invariants take in general values in a field of rational
functions in one or two formal variables and are conjecturally related to the quantum
Donaldson-Thomas invariants of Kontsevich and Soibelman by a refined
multicover formula.

In order to make contact with previous
results, note that refined wallcrossing formulas have been derived in physical theories
defined by a Seiberg-Witten curve in \cite{DM-crossing,DG,CV-I,DGS,laminations}, and conjectured to
hold in more general situations. Moreover, motivic wallcrossing formulas for
Donaldson-Thomas invariants of quivers with potential have been
 also announced in \cite{KS-review}.
The wallcrossing formulas conjectured in this paper for refined generalized
Donaldson-Thomas invariants, are related to those of \cite{DG,CV-I,DGS,laminations}
by a refined multicover formula, as explained in more detail below.
In addition, it is worth noting that
the invariants conjectured here are also equivariant residual invariants
with respect to a torus action. Therefore a rigorous construction
would require an equivariant localization theorem for motivic Donaldson-Thomas
invariants. Although the conjectures below are specifically
formulated for ADHM sheaves, analogous conjectures can be formulated  in more
general situations including abelian categories of coherent sheaves or coherent
perverse sheaves on Calabi-Yau threefolds. Previous results and conjectures in the mathematics literature are presented in \cite{motivic-degzero, N-refined}.

The main application of the conjectures in section
(\ref{motivicconjectures})  is a recursive formula presented
in section (\ref{recursionformula}). This formula determines
the Poincar\'e and Hodge polynomial of moduli spaces of Hitchin pairs with coprime
rank and degree in terms of asymptotic motivic ADHM invariants. The later are in turn
determined by string theoretic techniques, the results being summarized in
section (\ref{asympmotivicsect}). In section (\ref{computations}) it is checked by direct
computation that the
resulting expressions are in agreement with the results of \cite{hitchin-selfd,Bettinumbers,
HRV, Mirror-Hodge} in many concrete examples.
This provides strong evidence for the validity of the conjectural formalism proposed here.
Note that Higgs sheaves on curves are also employed in
\cite{BPS-Higgs} as a computational device for local BPS invariants of toric surfaces.

\subsection{Refined Wallcrossing Conjectures}\label{motivicconjectures}

In order to fix the notation, let $\Delta(r,e)\subset \IR_{>0}$ be the (finite)
set of positive critical stability parameters of type $(r,e)\in \IZ_{\geq 1}\times \IZ$.
For any $n\in \IZ$, and any formal variable $y$ let
\[
[n]_y = {y^n-y^{-n}\over y-y^{-1}}\in \IQ(y)
 \]
\begin{conj}\label{motivicinvariants}
Let $\gamma=(r,e)\in \IZ_{\geq 1}
\times \IZ$. Then
there exist refined equivariant residual ADHM invariants $A_\delta(r,e)(y)\in
\IQ(y)$, for any $\delta\in \IR$,
 and refined equivariant residual Higgs sheaf invariants $H(r,e)(y) \in \IQ(y)$
so that $A_\delta(r,e)(1) = A_\delta(r,e)$,
$H(r,e)(1)=H(r,e)$ and  the following wallcrossing formulas hold.

$(i)$ Let $\delta_c \in \Delta(r,e)$ be critical stability parameter and
$\delta_{c-}<\delta_c$,  $\delta_{c+}>\delta_c$ be noncritical
stability parameters so that $[\delta_{c-},\delta_c)\cap \Delta(r,e)=\emptyset$,
$(\delta_c,\delta_{c+}]\cap \Delta(r,e)=\emptyset$. The following
wallcrossing formula holds for $\delta_{c\pm}$ sufficiently close to $\delta_c$
\be\label{eq:wallB}
\bal
& A_{\delta_{c}+}(\gamma)(y) - A_{\delta_c-}(\gamma)(y) = \\
& \mathop{\sum_{l\geq 2}}_{} {1\over (l-1)!}
\mathop{\sum_{\gamma_1+\cdots+\gamma_l=\gamma}}_{\mu_{\delta_c}(\gamma_1)=\mu(\gamma_2)=\cdots = \mu(\gamma_l)}
A_{\delta_{c-}}(\gamma_1)\prod_{i=2}^l (-1)^{e_i-r(g-1)}[e_i-r_i(g-1)]_y
H(\gamma_i)(y)\\
\eal
\ee
where the sum in the right hand side of \eqref{eq:wallB} is finite.
Moreover
$[\delta_{c-},\delta_c)\cap \Delta(r_1,e_1)=\emptyset$,
$(\delta_c,\delta_{c+}]\cap \Delta(r_1,e_1)=\emptyset$
for all $\gamma_1=(r_1,e_1)$ in the right hand side of \eqref{eq:wallB}.

$(ii)$ Let
$\delta_{-}<0$,  $\delta_{+}>0$ be noncritical
stability parameters so that $[\delta_{-},0)\cap \Delta(r,e)=\emptyset$,
$(0,\delta_{+}]\cap \Delta(r,e)=\emptyset$. The following
wallcrossing formula holds for $\delta_{\pm}$ sufficiently close to $0$
\be\label{eq:wallA}
\bal
& A_{\delta_+}(\gamma)(y) - A_{\delta_-}(\gamma)(y) = \\
& \mathop{\sum_{l\geq 1}}_{} {1\over l!} \mathop{\sum_{\gamma_1+\cdots+\gamma_l=\gamma}}_{\mu(\gamma_i)=\mu(\gamma),\ 1\leq i \leq l} \prod_{i=1}^l
(-1)^{e_i-r(g-1)}[e_i-r_i(g-1)]_y H(\gamma_i)(y)\\
& +
\mathop{\sum_{l\geq 2}}_{} {1\over (l-1)!}
\mathop{\sum_{\gamma_1+\cdots+\gamma_l=\gamma}}_{\mu(\gamma_i)=\mu(\gamma),\ 1\leq i \leq l} A_{\delta_-}(\gamma_1)(y)\prod_{i=2}^l (-1)^{e_i-r(g-1)}
[e_i-r_i(g-1)]_y
H(\gamma_i)(y)\\
\eal
\ee
where the sum in the right hand side of \eqref{eq:wallA} is finite. Moreover,
$[\delta_{-},0)\cap \Delta(r_1,e_1)=\emptyset$,
$(0,\delta_{+}]\cap \Delta(r_1,e_1)=\emptyset$
for all $\gamma_1=(r_1,e_1)$ in the second line of the right hand side of
equation \eqref{eq:wallA}.

Moreover $A_\delta(r,e) \in \IZ[y,y^{-1}]$
if $\delta\in \IR$ is noncritical, and
$H(r,e)(y) \in \IZ[y,y^{-1}]$ if $(r,e)$ are coprime.
\end{conj}

As mentioned above the invariants $A_\delta(r,e)\in \IZ[y,y^{-1}]$, $H(r,e)(y)$
 are conjecturally related to residual equivariant Kontsevich-Soibelman invariants
${\overline A}_\delta(r,e)(y)\in \IZ[y,y^{-1}]$,
${\overline H}(r,e)(y)\in \IZ[y,y^{-1}]$
by a refined multicover formula. For $v=1$ invariants this formula states
simply that $A_\delta(r,e)(y)= {\overline A}_\delta(r,e)(y)$, while the explicit
formula for $v=0$ is given below.
\begin{conj}\label{motivicmulticover}
Under the same hypothesis as in conjecture (\ref{motivicinvariants}),
the following relation holds for any $(r,e)\in \IZ_{\geq 1} \times \IZ$
\be\label{eq:multicover}
H(r,e)(y) = \sum_{\substack{k\in \IZ,\ k\geq 1\\ k|r,\ k|e}} {1\over k\, [k]_y} {\overline H}
\bigg({r\over k}, {e\over k}\bigg)(y^k).
\ee

\end{conj}

The refined wallcrossing formulas \eqref{eq:wallB}, \eqref{eq:wallA}
are formal quantum generalizations of
the wallcrossing formulas derived in \cite{chamberII}. Refined, or quantum,
wallcrossing formulas have been physically derived
in \cite{DG,CV-I,DGS} using arguments analogous to \cite{DM-split}.
In particular a  refinement of the semiprimitive wallcrossing formula
of \cite{DM-split} has been formulated in \cite{DGS}.
A motivic wallcrossing formula  has been also  announced in \cite{KS-review}.
By analogy with \cite[Sect. 4]{chamberII}, \cite[Sect. 4]{higherrank},
the wallcrossing formulas conjectured in (\ref{motivicinvariants})
can be shown to agree with the refined semiprimitive wallcrossing formulas
of  \cite{CV-I,DGS,laminations},
once the multicover formula \eqref{eq:multicover} is properly taken into account.
In particular the above refined multicover formula can be easily inferred from
\cite[Sect 4.]{CV-I}.
The details are similar to those in \cite[Sect. 4]{chamberII}, \cite[Sect.  4]{higherrank},
hence will be omitted.

 Finally note that  a refined formula has been also derived in \cite{DM-crossing} for primitive wallcrossing using arguments analogous to \cite{DM-split},
and shown to be in a agreement with
wallcrossing formulas
for Poincar\'e and Hodge polynomials of moduli spaces of stable sheaves
on surfaces \cite{hodgenumbers,KY-Betti,chamber}.
The formula derived in \cite{DM-crossing} is in fact doubly refined,
the BPS states being simultaneously graded by spin and $U(1)_R$-charge
quantum numbers. This motivates the following further refinement of conjecture (\ref{motivicinvariants}), which can be physically justified using arguments analogous
to \cite{DM-split, DM-crossing,DGS}.

Let $(u,v)$ be formal variables, and $(u^{1/2}, v^{1/2})$ be formal square roots.
For any $n\in \IZ$ set
\[
 [n]_{(u,v)} = {(uv)^{n/2}-(uv)^{-n/2}\over (uv)^{1/2}-(uv)^{-1/2}}\in \IQ(u^{1/2},v^{1/2}).
 \]
\begin{conj}\label{bigradedinv}
Under the same conditions as in conjecture (\ref{motivicinvariants}) there exist
 doubly refined  equivariant residual ADHM invariants ${A}_\delta(r,e)(u,v)\in
\IQ(u^{1/2},v^{1/2})$,
 and doubly refined   Higgs sheaf invariants $H(r,e)(u,v) \in \IQ(u^{1/2},v^{1/2})$
so that

$(i)$ ${ A}_\delta(r,e)(u,u) = A_\delta(r,e)(u)$,
$H(r,e)(u,u)=H(r,e)(u)$,\\
${A}_\delta(r,e)(u,v)\in \IZ[u^{1/2},u^{-1/2}, v^{1/2}, v^{-1/2}]$
if $\delta$ is noncritical and $H(r,e)(u,v)\in
\IZ[u^{1/2},u^{-1/2}, v^{1/2}, v^{-1/2}]$
if $(r,e)$ are coprime.

$(ii)$ $A_\delta(r,e)(u,v)$ satisfy wallcrossing formulas obtained by
substituting \\
$A_\delta(\gamma_i)(u,v), H(\gamma_i)(u,v), [e_i-r_i(g-1)]_{(u,v)}$
for $A_\delta(\gamma_i)(y), H(\gamma_i)(y), [e_i-r_i(g-1)]_{y}$ in \eqref{eq:wallB}, \eqref{eq:wallA}.

$(iii)$ There exist alternative Higgs sheaf invariants ${\overline H}(r,e)(u,v)\in
\IZ[u^{1/2},u^{-1/2}, v^{1/2}, v^{-1/2}]$, $(r,e)\in \IZ_{\geq 1}\times \IZ$
so that ${H}(r,e)(u,v)$, ${\overline H}(r,e)(u,v)$ satisfy a multicover formula
obtained by making the same substitutions in \eqref{eq:multicover}.
\end{conj}

Note that the same notation $A_\delta(r,e)$, $H(r,e)$;
$A_\delta(r,e)(y)$,  $H(r,e)(y)$; $A_\delta(r,e)(u,v)$, $H(r,e)(u,v)$
is (abusively)  employed
for rational, respectively motivic and  refined motivic invariants.
By convention, the distinction will reside only in the number of arguments of
these rational functions. Therefore if no arguments are present, $A_\delta(r,e)$, $H(r,e)$
are rational numbers, if one argument is present they are rational functions of one variable
etc. Moreover, the invariants $H(r,e)(y)$ will be called refined Higgs invariants in the following. The invariants
$A_{\delta_\pm}(r,e)(y)$ with $\delta_\pm$ close to 0
as in (\ref{motivicinvariants}.$ii$)  will be denoted
by $A_{0\pm}(r,e)(y)$. Similarly the invariants $A_\delta(r,e)(y)$, with $\delta > \mathrm{max}\, \Delta(r,e)$
respectively $\delta < \mathrm{min}\, \Delta(r,e)$ will be denoted by $A_{\pm\infty}(r,e)(y)$ and
referred to as asymptotic invariants.

Finally note that the duality isomorphisms \eqref{eq:dualmoduliA}, \eqref{eq:dualmoduliB}
yield relations of the form
\be\label{eq:dualmotivic}
A_\delta(r,e)(y) =A_{-\delta}(r,-e+2r(g-1))(y)\qquad
 H(r,e)(y) =  H(r,-e+2r(g-1))(y)
 \ee
 for all $(r,e)\in \IZ_{\geq 1}\times \IZ$.
 Moreover, the isomorphisms \eqref{eq:shiftmoduli} imply that
 \be\label{eq:shiftmotivic}
 H(r,e)(y)=H(r,e+r)(y).
 \ee
for any $(r,e)\in \IZ_{\geq 1}\times \IZ$.
Therefore for fixed $r$ there are only $r$ a priori distinct invariants
$H(r,e)(y)$. Obviously entirely analogous formulas hold for the refined motivic
invariants $A_\delta(r,e)(u,v)$, $H(r,e)(u,v)$.

\subsection{Recursion formula for refined Higgs invariants}\label{recursionformula}
For the purpose of the present paper, the main application
 of conjectures  (\ref{motivicinvariants}), (\ref{bigradedinv})  is a
recursion formula for the invariants $H(r,e)(y)$, $H(r,e)(u,v)$ which
determines inductively all invariants $H(r,e)(y)$, $H(r,e)(u,v)$,
$(r,e)\in \IZ_{\geq 1}\times \IZ$
in terms of the asymptotic invariants $A_{+\infty}(r,e)(y)$,
$A_{+\infty}(r,e)(u,v)$.

 In the following $X$ is assumed to be a smooth projective curve of genus $g\geq 2$
 and $p=\mathrm{deg}(M_1)\geq 0$.
 For any $\gamma=(r,e)$, let $\wgamma =
(r, -e+2r(g-1))$, ${\widetilde e} = -e + 2r(g-1)$. For a stability parameter
$\delta$ let $\mu_\delta(\gamma)
=(e+\delta)/r$, $\mu(\gamma) = e/r$. Given $\gamma=\IZ\times \IZ$, the notation
$\gamma=(r(\gamma),e(\gamma))$ will also be used on occasion.

The recursion formula will be written in detail only for the refined  invariants
$H(r,e)(y)$ since the analogous formula for the doubly refined invariants
$H(r,e)(u,v)$ follows by obvious substitutions, as explained in conjecture (\ref{bigradedinv}).
 Let $\gamma=(r,e)\in \IZ_{\geq 1}\times \IZ$ be an arbitrary numerical type.
Then the following wallcrossing formula holds.
\be\label{eq:higgsrecursionA}
\bal
& (-1)^{e-r(g-1)} [e-r(g-1)]_y H(\gamma)(y) =A_{+\infty}({\gamma})(y)-A_{+\infty}(\wgamma)(y)\\
& +
\mathop{\sum_{l\geq 2}}_{}{(-1)^{l-1}\over (l-1)!}
\sum_{\substack{ \gamma_1,\ldots,\gamma_l \in \IZ_{\geq 1}\times \IZ\\
\gamma_1+\cdots+\gamma_l=\gamma\\
\mu(\gamma)< \mu(\gamma_i), \ 2\leq i\leq l\\}}
A_{+\infty}(\gamma_1)(y)
\prod_{i=2}^l (-1)^{e_i-r_i(g-1)}[e_i-r_i(g-1)]_y H(\gamma_i)(y)\\
& -\mathop{\sum_{l\geq 2}}_{}{(-1)^{l-1}\over (l-1)!}
\sum_{\substack{ \gamma_1,\ldots,\gamma_l \in \IZ_{\geq 1}\times \IZ\\
\gamma_1+\cdots+\gamma_l=\wgamma\\
 \mu(\wgamma)\leq \mu(\gamma_i), \ 2\leq i\leq l\\}}
A_{+\infty}(\gamma_1)(y)
\prod_{i=2}^l (-1)^{e_i-r_i(g-1)}[e_i-r_i(g-1)]_y H(\gamma_i)(y)\\
& -\mathop{\sum_{l\geq 2}}_{} {1\over l!}
\sum_{\substack{
\gamma_1,\ldots,\gamma_l \in \IZ_{\geq 1}\times \IZ\\
\gamma_1+\cdots+\gamma_l=\gamma\\
\mu(\gamma)=\mu(\gamma_i),\ 1\leq i\leq l\\}}
\prod_{i=1}^l (-1)^{e_i-r_i(g-1)}[e_i-r_i(g-1)]_y H(\gamma_i)(y)\\
\eal
\ee
where  the sum in the right hand side of equation
\eqref{eq:higgsrecursionA} contains only finitely many nontrivial terms.
The derivation of the recursion formula \eqref{eq:higgsrecursionA}
from the wallcrossing formulas \eqref{eq:wallB}, \eqref{eq:wallA} is
presented in section (\ref{recsection}).

\begin{rema}\label{recrema}
$(i)$ Note
that only invariants $H(r_i,e_i)(y)$ with $r_i<r$ enter the sum in right hand side
of \eqref{eq:higgsrecursionA}. Therefore this relation completely determines all
invariants $H(r,e)$, $(r,e)\in \IZ_{\geq 1}\times \IZ$ if all invariants $A_{+\infty}(r,e)(y)$ are known. A conjectural formula for the asymptotic refined  ADHM invariants
$A_{+\infty}(r,e)(y)$ will be derived in the next section using string duality.

$(ii)$ Given relations \eqref{eq:dualmotivic}, \eqref{eq:shiftmotivic},
equations \eqref{eq:higgsrecursionA} are in fact an overdetermined set of recursion
relations for refined  Higgs invariants. If conjecture (\ref{motivicinvariants}) holds,
all these equations are compatible, and one can choose the most economical one
for concrete computations.  In fact, one can obtain a simpler relation by taking
$e>2r(g-1)-c(r)$ in \eqref{eq:higgsrecursionA}. This results in $A_{+\infty}(\wgamma)=0$
and the second line in the right hand side is zero as well.
However, the simpler relation obtained this way is not necessarily the most efficient as far as computer time is concerned.
Concrete examples and computations
will be presented in section (\ref{computations}).
 \end{rema}

\subsection{Higgs invariants and cohomology of moduli spaces
of Hitchin pairs}\label{cohpairs}
The goal of this subsection is to formulate one more conjecture relating refined Higgs
invariants to the cohomology of moduli spaces of stable Hitchin pairs on $X$, for
coprime numerical invariants $(r,e)\in \IZ_{\geq }\times \IZ$. In the following
it is still assumed that the genus of $X$ is $g\geq 2$, and $p=\mathrm{deg}(M_1)
\geq 0$. Moreover,  $M_1\simeq \CO_X$ if $p=0$.

First recall that a Hitchin pair on $X$ with coefficient line bundle $L$
is a coherent sheaf $E$ equipped with a morphism $\Phi:E\to E\otimes_X L$.
The moduli theory of such objects has been extensively and intensively
studied in the mathematics literature
\cite{hitchin-selfd,MR1085642,infpairs,projective,simpsonI,simpsonII}. In particular, as recalled in section  (\ref{hitchinpairs}), there is a natural stability condition which yields
an algebraic  moduli stack ${\mathfrak H}(X,L,r,e)$ of finite type.
Moreover, suppose $\mathrm{deg}(L)\geq 2g-2$ and $L\simeq K_X$ if
$\mathrm{deg}(L) = 2g-2$.
There also exists a coarse moduli scheme $H^s(X,L,r,e)$
parameterizing isomorphism classes of stable objects.
If $(r,e)\in \IZ_{\geq 1}\times \IZ$ are coprime, any semistable Hitchin pair is stable,
and $H^s(X,L,r,e)$ will be denoted by $H(X,L,r,e)$.

The connection between Higgs sheaves and Hitchin pairs is based on
the observation that there is a natural forgetful morphism of moduli stacks
\[
{\mathfrak {Higgs}}(X,M_1,M_2,r,e) \to {\mathfrak{H}}(X,M_2^{-1},r,e)
\]
which simply forgets $\Phi_1:E\otimes_X M_1\to M_1$. Moreover, under the current assumptions, this morphism is compatible with stability for $(r,e)$ coprime,
and has a very simple structure as explained in section (\ref{hitchinpairs}). This
leads to the conjecture formulated below.

First note that for $(r,e)\in \IZ_{\geq 1}\times \IZ$ coprime, the degree of the
Poincar\'e polynomial $P_y(H(X,L,r,e))$ of the
smooth moduli space $H(X,L,r,e)$ is an even integer
$2m(r,e)$, $m(r,e)\in \IZ_{\geq 0}$. Under the same
conditions, let $H_{(u,v)}(H(X,L,r,e))$ denote the
Hodge polynomial of $H(X,L,r,e)$
(see \cite[Sect. 2.1]{HRV}, \cite[Sect. 2]{Mirror-Hodge} for definition
and properties.)
 \begin{conj}\label{Epolynomialpairs}
Under the above assumptions, let $L\simeq M_2^{-1}$.
Then
\be\label{eq:refEpolynomial}
\bal
 H(r,e)(y) & = (-1)^{e-r(g-1-p)} y^{-n(r,e)}P_{(-y)}(H(X,L,r,e))\\
  H(r,e)(u,v) & = (-1)^{e-r(g-1-p)} (uv)^{-n(r,e)/2}H_{(-u,-v)}(H(X,L,r,e))\\
  \eal
 \ee
 where  \[
 n(r,e) = r^2(g-1)+r(r-1)p + m(r,e).
 \]
   \end{conj}

\begin{rema}\label{Epolyrema}
$(i)$ The recursion relation \eqref{eq:higgsrecursionA} and conjecture (\ref{Epolynomialpairs}) determine all
 Hodge polynomials $ H_{(u,v)}(H(X,L,r,e))$
with $(r,e)\in \IZ_{\geq 1}\times \IZ$ coprime if the  asymptotic
refined  ADHM invariants
 are known for all $(r,e)\in \IZ_{\geq 1}\times \IZ$.  Conjectural formulas for
 these asymptotic invariants are
presented in the next subsection.

$(ii)$ Note that the recursion formula \eqref{eq:higgsrecursionA} determines in fact  all invariants $H(r,e)(y)$, $H(r,e)(u,v)$, including non-coprime pairs.
A priori, the Higgs invariants $H(r,e)(y)$ with $(r,e)$ not coprime are not related
in any direct way to the cohomology of moduli spaces of semistable Hitchin pairs
with the same numerical invariants. However, a conjectural relation based on the
multicover formula \eqref{eq:multicover} will be formulated  in the next subsection.
\end{rema}

\subsection{Asymptotic refined ADHM invariants}\label{asympmotivicsect}
As explained above, the invariants $H(r,e)(y)$, $H(r,e)(u,v)$
are completely determined by the recursion relation \eqref{eq:higgsrecursionA}
if all  asymptotic refined ADHM invariants are known. A conjectural
formula for the generating function of  asymptotic refined ADHM invariants
is derived from string duality in section (\ref{gaugesect}).  Basically, this
 generating function  is determined
by the Nekrasov partition function \cite{Nekrasov:2002qd} of a five dimensional
supersymmetric gauge theory.

As shown below, the resulting formula involves an infinite formal sum over
Young tableaus $Y$. In order to fix conventions,  note that a nonempty Young tableau
$Y$ is identified with a partition
\[
|Y| = Y_{1} + \cdots + Y_{l(Y)}
\]
where $Y$ denotes the total number of boxes of $Y$ and $l(Y)$ denotes the
number of rows. For any $1\leq i \leq l(Y)$, $Y_i$ denotes the length of the $i$-th
row, and $Y_1\geq Y_2\geq \cdots \geq Y_{l(Y)}$. Boxes of $Y$ will be labeled
by $(i,j)\in \IZ\times\IZ$, $1\leq i\leq l(Y)$, $1\leq j \leq Y_i$.

\begin{conj}\label{refasympinvI}
Let $\CX=(X,M_1,M_2)$ be a triple as above  and let
$p=\mathrm{deg}(M_1)$. Let
\be\label{asympmotivic}
\CZ_{+\infty}(\CX,r; \lambda,y) = \mathop{\sum_{e\in \IZ}}_{} \lambda^e
A_{+\infty}(r,e)(y)
\ee
be the generating function for the rank $r\in \IZ_{\geq 1}$
asymptotic refined ADHM invariants conjectured in
(\ref{motivicinvariants}).
Then
\be\label{eq:refC}
Z_{+\infty}(\CX,r;\lambda, y) =
 \mathop{\sum_{|Y|=r}}_{}
\Omega^{(g,p)}_{Y}(\lambda,y)
\ee
where
\be\label{eq:omegaformulaI}
\bal
\Omega^{(g,p)}_Y(\lambda,y)=  & (-1)^{p|Y|}y^{-p\sum_{(i,j)\in Y}(i+j-2)+(g-1)\sum_{(i,j)\in Y}(-2i+2j+1-2Y_i+Y^t_j)}\\
&\lambda^{-p\sum_{(i,j)\in Y}(-i+j)+(g-1)\sum_{(i,j)\in Y}(2i+2j-1-2Y_i-Y^t_j)}\\
&\prod_{(i,j)\in
Y}F(\lambda^{-i-j+Y_i+Y^t_j+1}y^{i-j+Y_i-Y^t_j},y)\eal
\ee
and
\[
F(q,z)=
z^{1-g}\frac{(1-q)^{2g}}{(1-qz)(1-qz^{-1})}.
\]
By convention $\Omega^{(p)}_{\emptyset}(\lambda,y)=1$.
\end{conj}

The generating function of asymptotic doubly refined ADHM invariants
\be\label{asympmotivicII}
\CZ_{+\infty}(\CX,r; \lambda,u,v) = \mathop{\sum_{e\in \IZ}}_{} \lambda^e
A_{+\infty}(r,e)(u,v)
\ee
is conjecturally  determined as follows.
\begin{conj}\label{refasympinvII}
Under the same conditions as in conjecture (\ref{refasympinvI}),
\be\label{eq:refE}
Z_{+\infty}(\CX,r;\lambda, u,v) =
 \mathop{\sum_{|Y|=r}}_{}
\Omega^{(g,p)}_{Y}(\lambda,u,v)
\ee
where
\be\label{eq:omegaformulaII}
\bal
\Omega^{(g,p)}_Y(\lambda,u,v)= &(-1)^{p|Y|}(uv)^{-p\sum_{(i,j)\in Y}(i+j-2)/2+(g-1)\sum_{(i,j)\in Y}(-2i+2j+1-2Y_i+Y^t_j)/2}\\
&\lambda^{-p\sum_{(i,j)\in Y}(-i+j)+(g-1)\sum_{(i,j)\in Y}(2i+2j-1-2Y_i-Y^t_j)}\\
&\prod_{(i,j)\in
Y}G(\lambda^{-i-j+Y_i+Y^t_j+1}(uv)^{(i-j+Y_i-Y^t_j)/2},(uv)^{1/2},(uv^{-1})^{1/2})\eal
\ee
and
\[
G(q,z,w)=
z^{(1-g)}\frac{(1-qw)^{g}(1-qw^{-1})^g}{(1-qz)(1-qz^{-1})}.
\]
By convention $\Omega^{(p)}_{\emptyset}(\lambda,u,v)=1$.
\end{conj}

Using the recursion relation \eqref{eq:higgsrecursionA} and conjectures
(\ref{Epolynomialpairs}), (\ref{refasympinvI}), (\ref{refasympinvII}),
one can derive explicit formulas for the Hodge polynomials of the
moduli spaces $H(X,L,r,e)$ with $(r,e)$-coprime. Note in particular that
formulas \eqref{eq:refC}, \eqref{eq:refE} imply that all invariants $A_{+\infty}(\gamma_1)(y)$
in the right hand side of equation \eqref{eq:higgsrecursionA} are trivial if
$\mu(\gamma_1) < -(r-1)(2g-2+p)$.
Concrete computations
are presented in section (\ref{computations}) for $r=1,2,3$
and various values of $g\geq 2$, $p\geq 0$. In all cases, the resulting formulas
are in agreement with the direct localization computations of Hitchin \cite{hitchin-selfd},
Gothen \cite{Bettinumbers} as well as the Hausel-Rodriguez-Villegas formula
\cite{HRV,Mirror-Hodge}.  A brief survey of the results in the  mathematics
literature on the subject is presented in appendix (\ref{survey}).
Moreover, direct computations in all examples considered in section (\ref{computations})
support the following intriguing conjecture.

\begin{conj}\label{intmotivic}
Under the same conditions as in conjecture (\ref{refasympinvI}),
for fixed $r\geq 1$,
the refined invariants ${\overline H}(r,e)(y)$, ${\overline H}(r,e)(u,v)$
are independent of $e\in \IZ$. In particular, they take the same value for all pairs $(r,e)$,
coprime or not.
\end{conj}
In fact, since the first version of this work was posted, the recursion relation 
\eqref{eq:higgsrecursionA} has been beautifully solved by Mozgovoy in 
\cite{adhm-rec}, and the solution has been proven to be in agreement with the Hausel-Rodriguez-Villegas 
invariants. Furthermore, Mozgovoy's solution also satisfies the multicover 
formula \eqref{eq:multicover} and has the property stated in conjecture (\ref{intmotivic}).

{\it Acknowledgments.} We are very grateful to Ugo Bruzzo,
Ron Donagi, Daniel Jafferis, Yunfeng Jiang, Dominic Joyce,
Greg Moore, Artan Sheshmani, Balasz Szendroi, Chris Woodward,
and especially Tamas Hausel, Ludmil Katzarkov, Sergey Mozgovoy,
Tony Pantev and Fernando Rodriguez-Villegas
for  their interest in this work and many helpful discussions. 
We owe special thanks to Sergey Mozgovoy for sending us 
his paper \cite{adhm-rec} before publication. 
DED would also like to thank the
organizers of VBAC 2009 Berlin for an excellent mathematical atmosphere which
prompted the research reported here. The work of
WYC is supported by DOE grant DE-FG02-96ER40959.  The work of DED
was partially supported by NSF grant PHY-0854757-2009.

\newpage
\section{ADHM invariants, Hitchin pairs and wallcrossing}

\subsection{Review of ADHM sheaves}\label{review}
Let $X$ be a smooth projective curve over $\IC$ of genus $g\geq 2$.
Let $M_1,M_2$ be line bundles on $X$ so that $M_1\otimes_X M_2\simeq K_X^{-1}$, and fix such an isomorphism in the following. Let $\mathrm{deg}(M_1)=p$,
$\mathrm{deg}(M_2) = -2g-2-p$, $p\in \IZ$ and $\CX=(X,M_1,M_2)$.

The abelian category $\CC_\CX$ of ADHM sheaves was defined in
\cite[Sect. 3]{chamberI}
as follows. The objects of $\CC_\CX$ are collections $\CE=(E, V,
\Phi_1,\Phi_2,\phi,\psi)$ on $X$ where $E$ is a coherent sheaf on $X$,
$V$ is a finite dimensional complex vector space, and
$\Phi_{i}:E\otimes_X M_i \to E$, $i=1,2$ ,
$\phi:E\otimes_X M_1\otimes_X M_2 \to V\otimes \CO_X$, $
\psi:V\otimes \CO_X\to E$ are morphisms of $\CO_X$-modules
satisfying the ADHM relation
\be\label{eq:ADHMrelation}
\Phi_1\circ(\Phi_2\otimes 1_{M_1}) - \Phi_2\circ (\Phi_1\otimes 1_{M_2})
+ \psi\circ \phi =0.
\ee
The morphisms of $\CC_\CX$ are natural morphisms of quiver sheaves.

An  object $\CE$ of $\CC_{\CX}$ will be called locally free if $E$ is a coherent
locally free $\CO_X$-module.
Given a coherent $\CO_X$-module $E$  we will denote by $r(E)$, $d(E)$,
$\mu(E)$ the rank, degree, respectively slope of $E$ if $r(E)\neq 0$.
The type of an object $\CE$ of $\CC_{\CX}$ is the collection
$(r(\CE), d(\CE), v(\CE))= (r(E),d(E),\mathrm{dim}(V)))\in
\IZ_{\geq 0}\times \IZ\times \IZ_{\geq 0}$.

Note that the objects of $\CC_{\CX}$ with $v(\CE)=0$ are triples
$\CE=(E,\Phi_1,\Phi_2)$ so that
\be\label{eq:Higgsrelation}
\Phi_1\circ(\Phi_2\otimes 1_{M_1}) - \Phi_2\circ (\Phi_1\otimes 1_{M_2}) =0.
\ee
and form a full abelian subcategory of $\CC_\CX$.
These are known as Higgs sheaves on $X$ with coefficient bundle $M_1\oplus M_2$
(see \cite[App. A]{chamberI} for a brief summary of definitions and properties.)

The dual of a locally free ADHM sheaf $\CE=(E,V,\Phi_1,\Phi_2,
\phi,\psi)$ is defined by
\be\label{eq:dualADHM}
\begin{aligned}
{\widetilde E} & = E^\vee \otimes_X M^{-1}\\
{\widetilde \Phi}_i & = (\Phi_i^\vee\otimes 1_{M_i}) \otimes 1_{M^{-1}}: {\widetilde E}
\otimes
M_i \to {\widetilde E} \\
{\widetilde \phi} & = \psi^\vee \otimes 1_{M^{-1}} : {\widetilde E}\otimes_X {M} \to
V^\vee\otimes \CO_X
\\
{\widetilde \psi} & = \phi^\vee : V^\vee \otimes \CO_X \to {\widetilde E} \\
\end{aligned}
\ee
where $i=1,2$.
Obviously, if $\CE$ is of type $(r,e,v)$,  ${\widetilde \CE}$ is of type $(r,-e+2r(g-1),v)$.

Any real parameter $\delta \in \IR$ determines a stability condition on $\CC_\CX$
\cite{chamberI,ranktwo}.
An object $\CE$ of $\CC_\CX$ is $\delta$-(semi)stable if any proper nontrivial
subobject $0\subset \CE'\subset \CE$ satisfies the inequality
\be\label{eq:deltastabA}
r(\CE) ( d(\CE') + \delta v(\CE'))\ (\leq)\  r(\CE') ( d(\CE) + \delta v(\CE)).
\ee
Standard arguments show that the $\delta$-stability condition satisfies the
 Harder-Narasimhan as well as Jordan-H\"older property for any $\delta\in \IR$.
Moreover the following properties hold for any  object
$\CE=(E,V,\Phi_1,\Phi_2\phi,\psi)$
of $\CC_\CX$ with $r(\CE)\geq 1$ and $v(\CE)=1$ \cite[Sect 3]{chamberI}
\begin{itemize}
\item[$(S.1)$] If $\CE$ is $\delta$-semistable for some $\delta\in \IR$, then
$\CE$ is locally free. In addition, if $\delta>0$ then $\psi$ is not identically zero;
if $\delta<0$, $\phi$ is not identically zero.
\item[$(S.2)$] If $\CE$ is $\delta$-stable for some $\delta \in \IR$,
the endomorphism ring of $\CE$ in $\CC_\CX$ is canonically isomorphic to $\IC$.
\item[$(S.3)$] $\CE$ is $\delta$-(semi)stable if and only if the dual ${\widetilde \CE}$
is $(-\delta)$-(semi)stable.
\end{itemize}
One also has the following boundedness results
 \cite[Lemm. 2.6, Lemm. 2.7, Cor. 2.8]{chamberI}
\begin{itemize}
\item[$(B.1)$] The set of isomorphism classes of locally free ADHM sheaves
of fixed type $(r,e,1)$ which are $\delta$-semistable for some $\delta\in \IR$
is bounded.
\item[$(B.2)$] For any $r\geq 1$ there exists an integer $c(r)\in \IZ$ so that
any $\delta$-semistable ADHM sheaf of type $(r,e,1)$ for some $\delta>0$
satisfies $e\geq c(r)$. Note that the integer $c(r)$  is not unique
unless required to be optimal with this property.
In fact the proof of \cite[Lemm. 2.6]{chamberI}
implies that any integer $$c(r) \leq -2(r-1)^2\mathrm{max}\{|\mathrm{deg}(M_1)|,
|\mathrm{deg}(M_2)|\}$$ satisfies this condition.
\end{itemize}
Note that for $v=0$ objects, $\delta$-stability is independent of $\delta$
and reduces to standard slope stability for Higgs sheaves on $X$.

A straightforward corollary of the above results is the existence of an algebraic
moduli stack of finite type $\mfm_\delta^{ss}(\CX,r,e)$ of $\delta$-semistable ADHM
sheaves on $X$ of type $(r,e,1)$ for any $(r,e)\in \IZ_{\geq 1}\times \IZ$
and any $\delta\in \IR$. The substack $\mfm_\delta^s(\CX,r,e)$
of $\delta$-stable objects is separated
and has the structure of a $\IC^\times$-gerbe over an algebraic moduli space
$M_\delta^{ss}(\CX,r,e)$.
Property $(S.3)$ also yields a canonical isomorphism
\be\label{eq:dualmoduliA}
\mfm^{ss}_{\delta}(\CX,r,e)\simeq \mfm^{ss}_{\delta}(\CX,r,-e+2r(g-1))
\ee
for any $(r,e)\in \IZ_{\geq 1}\times \IZ$
and any $\delta\in \IR$.

Moreover there is a stability chamber structure on $\IR_{>0}$
as follows \cite[Sect. 4]{chamberI}.
For a fixed type $(r,e)\in \IZ_{\geq 1}\times \IZ$, three exists a finite
set $\Delta(r,e) \subset \IR_{>0}$ of critical stability parameters so that
\begin{itemize}
\item[$(C.1)$] For any $\delta \in \IR_{>0}\setminus \Delta(r,e)$,
$\delta$-semistability is equivalent to $\delta$-stability i.e.
$\mfm^{ss}_{\delta}(\CX,r,e)=\mfm^{s}_{\delta}(\CX,r,e)$.
\item[$(C.2)$] For any $\delta > \rm{max}\, \Delta(r,e)$ $\delta$-stability is equivalent
with the following asymptotic stability condition. An object
$\CE=(E,V,\Phi_i, \phi,\psi)$ with $v=1$
is asymptotically stable if $E$ is locally free, $\psi$ nontrivial, and there is no
proper saturated subsheaf $0\subset E'\subset E$ preserved by $\Phi_i$, $i=1,2$
so that $Im(\psi)\subseteq E'$.
\end{itemize}

Finally note that there is a torus ${\bf S}=\IC^\times$ action on the moduli stacks
$\mfm^{ss}_{\delta}(\CX,r,e)$ so that
\be\label{eq:toractA}
t\times (E,V,\Phi_1,\Phi_2, \phi,\psi)\to (E,V,t^{-1}\Phi_1,t^{}\Phi_2,\phi,\psi)
\ee
on closed points. According to \cite[Thm. 1.5]{chamberI}, for noncritical stability
parameter $\delta\in \IR_{>0}\setminus \Delta(r,e)$, the stack theoretic fixed locus
$\mfm_{\delta}^{ss}(\CX,r,e)^{\bf S}$
is universally closed over $\IC$.
Moreover, the algebraic moduli space $M_{\delta}^{ss}(\CX,r,e)$
has a perfect obstruction theory. Therefore residual $\delta$-ADHM invariants
$A_{\delta}(r,e)\in \IZ$ can be  defined in each chamber by equivariant virtual localization.
Wallcrossing formulas for these invariants have been derived in \cite[Thm. 1.1]{chamberII}
using Joyce-Song theory \cite{genDTI}.

For future reference note that there is a completely
analogous torus action on the moduli stack $\obj(\CC_\CX)$
 of all objects of $\CC_\CX$, which
is an algebraic stack of locally finite type over $\IC$.
In particular, this yields a torus action on the moduli stack
${\mathfrak {Higgs}}^{ss}(\CX,r,e)$ of slope-semistable
Higgs sheaves on $X$ with fixed $(r,e)\in \IZ_{\geq 1}\times \IZ$,
which is an algebraic stack of finite type over $\IC$. The wallcrossing formulas
in \cite[Thm. 1.1]{chamberII} are written in terms residual equivariant generalized Donaldson-Thomas invariants $H(r,e)\in \IQ$ defined via Joyce-Song theory
applied to the stacks ${\mathfrak {Higgs}}^{ss}(\CX,r,e)$. For curves $X$ of
genus $g\geq 1$, the invariants $H^{\bf S}(r,e)$ are trivial, hence the wallcrossing
formulas state that the invariants $A_\delta(r,e)$ are independent of
$\delta$.

In order to conclude this section, note that the stacks ${\mathfrak {Higgs}}^{ss}(\CX,r,e)$
have the following simple properties.
By analogy with \eqref{eq:dualmoduliA}, there is  a canonical torus equivariant isomorphism
\be\label{eq:dualmoduliB}
{\mathfrak {Higgs}}^{ss}(\CX,r,e) \simeq {\mathfrak {Higgs}}^{ss}(\CX,r,-e+2r(g-1))
\ee
In addition, taking tensor product by a fixed degree one line bundle on $X$ yields
an equivariant isomorphism
\be\label{eq:shiftmoduli}
{\mathfrak {Higgs}}^{ss}(\CX,r,e) \simeq {\mathfrak {Higgs}}^{ss}(\CX,r,e+r)
\ee
for any $(r,e)\in \IZ_{\geq 1}\times \IZ$. Finally note that for $(r,e)$ coprime slope
semistability is equivalent to slope stability, and
the stack ${\mathfrak {Higgs}}^{ss}(\CX,r,e)$ has a $\IC^\times$-gerbe structure
over a quasi-projective scheme $Higgs^{ss}(\CX,r,e)$.

\subsection{Connection with Hitchin pairs}\label{hitchinpairs}
Let $L$ be a fixed line bundle on $X$.
Recall that a Hitchin pair \cite{hitchin-selfd,MR1085642}
on $X$ with coefficient bundle $L$ is defined is a pair $(E,\Phi)$ where
$E$ is a coherent sheaf on $X$ and $\Phi:E\to E\otimes_X L$ a morphism
of coherent sheaves. Such a pair is called (semi)stable if any proper nontrivial
subsheaf $0\subset E'\subset E$ so that $\Phi(E') \subset E'\otimes_X L$
satisfies the inequality
\be\label{eq:HitchinstabA}
r(E) d(E') \ (\leq) \ r(E') d(E).
\ee
Note that if $r(E)>0$, semistability implies that $E$ is locally free.
In the following $L$ be either $K_X$ or a line bundle on $X$ of degree
$d(L)>2g-2$. This will be implicitly assumed in all statements below.

Well-known results in the literature
\cite{hitchin-selfd,MR1085642,infpairs,projective,simpsonI,simpsonII} establish the existence of an algebraic
stack of finite type ${\mathfrak H}(X,L,r,e)$ of semistable Hitchin pairs
of fixed type $(r(E),d(E))=(r,e)\in \IZ_{\geq 1}\times \IZ$. Moreover, if $(r,e)$ are
coprime, this stack is a $\IC^\times$-gerbe over a smooth quasi-projective
variety $H(X,L,r,e)$. For $L=K_X$, $H(X,L,r,e)$ is commonly referred to as the
Hitchin integrable system.

Note that there is a torus $\IC^\times$ action on the stack  ${\mathfrak H}(X,L,r,e)$
 given by $t\times (E,\Phi)\to (E,t^{-1}\Phi)$ on closed points.
The stack theoretic fixed locus is universally closed. In particular, for $(r,e)$ coprime,
there is an induced torus action on the moduli scheme $H(X,L,r,e)$, and the
fixed locus is a smooth projective scheme over $\IC$.

The relation between ADHM sheaves and Hitchin pairs is summarized in the following
simple observations.
\begin{itemize}
\item[$(AH.1)$] Suppose $M_1=\CO_X$, $M_2=K_X^{-1}$ and let $(r,e)\in \IZ_{\geq 1}
\times \IZ$ be coprime. Then there is an isomorphism
\be\label{eq:adhmhitchinA}
{\mathfrak {Higgs}}(\CX,r,e)\simeq \IC \times {\mathfrak{H}}(X,K_X,r,e).
\ee
\item[$(AH.2)$] Suppose $M_2$ is a line bundle of degree $2-2g-p$, where
$p\in \IZ_{>0}$. Then there is an isomorphism
\be\label{eq:adhmhitchinB}
{\mathfrak {Higgs}}(\CX,r,e)\simeq  {\mathfrak{H}}(X,M_2^{-1},r,e).
\ee
\end{itemize}

Both statements rely on the fact that for coprime $(r,e)$ slope semistability is
equivalent to slope stability. Therefore the endomorphism ring of any
semistable object $\CE$ is canonically isomorphic to $\IC$.

Then note that in the first case,
given any semistable object $\CE=(E,\Phi_1,\Phi_2)$ the relation \eqref{eq:Higgsrelation}
implies that $\Phi_1:E\to E$ is an endomorphism of $\CE$ since it obviously commutes
with itself. Therefore it must be of the form $\Phi_1= \lambda 1_E$ for some
$\lambda \in \IC$. In particular, it preserves any subsheaf $E'\subset E$. Generalizing this
observation to
flat families it follows that there is an forgetful morphism
\[
{\mathfrak {Higgs}}(\CX,r,e)\to {\mathfrak{H}}(X,K_X,r,e)
\]
projecting $(E,\Phi_1,\Phi_2)$ to $(E, \Phi_2\otimes 1_{K_X})$. The isomorphism
\eqref{eq:adhmhitchinA} then follows easily.

In the second case, note that given a semistable Higgs sheaf $(E,\Phi_1,\Phi_2)$,
of type $(r,e)$, the data
\[
\CE'=\left(E\otimes_X M_1^{-1}, \Phi_1\otimes 1_{M_1^{-1}}, \Phi_2\otimes 1_{M_1^{-1}}
\right)
\]
determines a semistable Higgs sheaf  of type $(r,e-r\mathrm{deg}(M_1))=(r,e-rp)$.
Relation \eqref{eq:Higgsrelation} implies that $\Phi_1\otimes 1_{M_1^{-1}}$
is a morphism of (semistable) Higgs sheaves. However $\mu(\CE)>\mu(\CE')$
since $p>0$, therefore any such morphism must vanish. This completes the proof.

\subsection{Remarks on refined wallcrossing conjectures}\label{remarksect}
This subsection consists of several remarks on  conjectures (\ref{motivicinvariants})
(\ref{bigradedinv}). It can be skipped with no loss of essential information.

$(i)$ First note that given any two objects $\CE_1,\CE_2$ of $\CC_\CX$ with
$v(\CE_1)+v(\CE_2)\leq 1$, it has been proven in \cite[Lemm. 7.4]{chamberI}
that the expression
\be\label{eq:altsum}
\mathrm{dim}\, \mathrm{Ext}^0_{\CC_\CX}(\CE_1,\CE_2)
-\mathrm{dim}\, \mathrm{Ext}^1_{\CC_\CX}(\CE_1,\CE_2) -
\mathrm{dim}\, \mathrm{Ext}^0_{\CC_\CX}(\CE_2,\CE_1)
+ \mathrm{dim}\, \mathrm{Ext}^1_{\CC_\CX}(\CE_2,\CE_1)
\ee
depends only on the numerical types of the two objects.
Moreover, if $\CE_1,\CE_2$ determine closed points in the stack theoretic fixed locus
$\obj(\CC_\CX)^{\bf S}$, there is an induced torus action on all the extension
groups in \eqref{eq:altsum} and the same statement holds for the alternating sum
of dimensions of fixed, respectively moving parts.
This technical condition makes both Joyce-Song and
Kontsevich-Soibelman theories applicable to non-Calabi-Yau categories, which is the
present case.

$(iii)$ As pointed out in \cite{DG}, the quantum Donaldson-Thomas invariants of
Kontsevich and Soibelman can be naturally identified with the refined topological
string invariants constructed in \cite{IKV} via the refined topological vertex formalism.
The asymptotic  invariants $A_{\pm \infty}(r,e)(y)$ are refinements
of the integral invariants $A_{\pm \infty}(r,e)$, which are in turn identical to
local stable pair invariants according to \cite{modADHM}. Therefore it entirely natural
to expect these invariants to be determined by the refined BPS counting invariants
of a local curve. The later can be inferred from the Nekrasov partition function
of a five dimensional gauge theory as explained in section (\ref{gaugesect}).

$(v)$ Finally note that assuming an equivariant localization result for motivic
invariants one can conjecture more refined wallcrossing formulas
for the residual contributions of individual components of the fixed loci. This follows
from the stack function relations derived in \cite[Sect. 3]{chamberII}.

\subsection{Derivation of recursion formula}\label{recsection}
The purpose of this section is to prove the recursion relation \eqref{eq:higgsrecursionA},
given the wallcrossing formulas \eqref{eq:wallB}, \eqref{eq:wallA}.
The proof  is analogous to the proof of \cite[Lemm. 3.8]{chamberII}.
The main steps will be outlined below for completeness.

According to property $(B.2)$ in section (\ref{review})
for any fixed $r\geq 1$ there
exists an integer $c(r)\in \IZ$ so that
all invariants $A_\delta(r,e)(y)$, for any $\delta>0$,
are identically zero if $e<c(r)$. Moreover, this integer is not unique unless
required to be optimal with this property; any integer $c(r)\leq -(r-1)^2(2g-2+p)$
satisfies this condition.  In the following set
\be\label{eq:lowerdegbound}
c(r) = -r(r-1)(2g-2+p) \qquad c(r') = -r'(r-1)(2g-2+p)
\ee
for any $r\in \IZ_{\ge q1}$,
$1\leq r'\leq r$. This is not an optimal choice, but it will facilitate the derivation of
formula \eqref{eq:higgsrecursionA}, as shown below.

Next note that
the wallcrossing formula \eqref{eq:wallB} is equivalent to
\be\label{eq:wallBrev}
\bal
& A_{\delta_{c}-}(\gamma)(y) - A_{\delta_c+}(\gamma)(y) = \\
& \mathop{\sum_{l\geq 2}}_{} {(-1)^{l-1}\over (l-1)!}
\mathop{\sum_{\gamma_1+\cdots+\gamma_l=\gamma}}_{\mu_{\delta_c}(\gamma_1)=\mu(\gamma_2)=\cdots = \mu(\gamma_l)}
A_{\delta_{c+}}(\gamma_1)\prod_{i=2}^l (-1)^{e_i-r(g-1)}[e_i-r_i(g-1)]_y
H(\gamma_i)(y)\\
\eal
\ee
 For any $n\in \IZ_{\geq 1}$ and any collection of $n$ positive integers
 $(l_1, \ldots, l_n) \in \IZ_{\geq 1}^{n}$,
define
 \be\label{eq:zerotoinftysetB}
\begin{aligned}
& {\sf S}^{(l_1,\ldots, l_n)}_{0,+\infty}(\gamma) =
\bigg\{ (\gamma_1, \eta_{1,1}, \ldots, \eta_{1,l_1},
 \ldots, {\eta}_{n,1}, \ldots, {\eta}_{n,l_n})\in
 (\IZ_{\geq 1}\times \IZ)^{\times( l_1+\ldots + l_n +1)}\,
 \bigg|\\
& {\gamma}_1+\sum_{i=1}^n\sum_{j=1}^{l_i} {\eta}_{i, j} =\gamma,\  \
\mu_0(r) \leq \mu(\gamma) <  \mu(\eta_{1,1})= \cdots = \mu(\eta_{1,l_1})< \\
& \mu(\eta_{2,1})= \cdots = \mu(\eta_{2,l_2})< \cdots <
\mu({\eta}_{n,1})= \cdots = \mu({\eta}_{n,l_n})< \mu_\delta(\gamma), \
\mu_0(r)\leq \mu(\gamma_1)
 \bigg\}\\
\end{aligned}
\ee
where $\mu_0(r)=c(r)/r$.
Then it straightforward to check that
the union
\be\label{eq:decompunion}
\bigcup_{n\geq 1}\bigcup_{l_1,\ldots, l_n \geq 1} {\sf S}^{(l_1,\ldots, l_n)}_{0,+\infty}(\gamma)
\ee
is a finite set.

Let $({\gamma}_{1}, {\eta}_{1,1}, \ldots, {\eta}_{1,l_1},
 \ldots, {\eta}_{n,1}, \ldots, {\eta}_{n,l_n}) \in
 {{\sf S}}^{(l_1,\ldots, l_n)}_{+,\delta}(\gamma)$
 be an arbitrary element, for some $n\geq 1$
 and $l_1,\ldots, l_n\geq 1$.
  Let $\mu_i$, $1\leq i\leq n$ denote the common value of the slopes
  $\mu({\eta}_{i,j})$,
 $1\leq j\leq l_i$.
 If $n\geq 2$, let also
 \[
 {\gamma}_{n-i+2} = {\gamma}_1 + {\eta}_{i,1} + \cdots + {\eta}_{n, l_{n}}
 \]
for $2\leq i\leq n$.
  Define the stability parameters $\delta_i$, $1\leq i\leq n$ by
 \be\label{eq:crtseqA}
 \begin{aligned}
 \mu_{\delta_1}({\gamma}_1)&  = \mu_n\\
 \mu_{\delta_i}(\gamma_i) & = \mu_{n+1-i}, \qquad
 2\leq i\leq n \qquad (\mathrm{if}\ n\geq 2). \\
 \end{aligned}
 \ee
 By construction, $\delta_{i}$ is a critical stability parameter of type $\gamma_{i}$
 for all $1\leq i\leq n$.
  Given the slope inequalities in \eqref{eq:zerotoinftysetB},
  it is straightforward to check that
  \be\label{eq:decreasingseq}
  0<\delta_n< \delta_{n-1} <\cdots < \delta_1.
  \ee
  Moreover, $\mu(\gamma_i)\geq \mu_0(r)$ for all $1\leq i\leq n$
  since
  the integers $c(r')$, $1\leq r'\leq r$ defined in \eqref{eq:lowerdegbound}
satisfiy
\be\label{eq:commslopeA}
{c(r')\over r'} = -(r-1)(2g-2+p)=\mu_0(r).
\ee

  Next note that the set $\Delta_\gamma$
  of all stability parameters constructed this way, for all $n\geq 1$ and
  any possible values of $l_1,\ldots, l_n$ is finite, since the set \eqref{eq:decompunion}
  is finite. Therefore one can choose stability parameters
  $0<\delta_{0+}< \mathrm{min}\,
  \Delta_\gamma$, $\delta_{+\infty}> \mathrm{max}\, \Delta_\gamma$.
  By construction $\Delta_\gamma$ contains all possible decreasing finite sequences of
  stability parameters of the form \eqref{eq:decreasingseq} with the property that
  there exists
  $$({\gamma}_{1}, {\eta}_{1,1}, \ldots, {\eta}_{1,l_1},
 \ldots, {\eta}_{n,1}, \ldots, {\eta}_{n,l_n}) \in
 (\IZ_{\geq }\times \IZ)^{\times(l_1+\cdots +l_n+1)}$$
 for some $l_1,\ldots, l_n\geq 1$ so that
 \begin{itemize}
\item[$(a)$] ${\gamma}_1+{\eta}_{1,1}+\cdots +{\eta}_{n,l_n}=\gamma$
\item[$(b)$] Conditions \eqref{eq:crtseqA} hold.
\end{itemize}
   In conclusion, successive applications of the wallcrossing formula \eqref{eq:wallBrev}
   yield
\be\label{eq:wallBB}
\bal
& A_{0+}(\gamma) - A_{+\infty}(\gamma) = \\
& \mathop{\sum_{n=1}^\infty}_{} \mathop{\sum_{l_1,\ldots,l_n\geq 1}}_{}
\prod_{i=1}^n {(-1)^{l_i}\over l_i!}
\sum_{\substack{
\gamma_1+\eta_{1,1}+\cdots + \eta_{1,l_1} + \cdots + \eta_{n,1}+
\cdots +\eta_{n,l_n}=\gamma,\\
\mu_0(r) \leq \mu(\gamma) < \mu(\eta_{1,1})=\cdots \mu(\eta_{1,l_1})<\cdots < \mu(\eta_{n,1})=\cdots =\mu(\eta_{n,l_n})\\
\mu_0(r)\leq \mu(\gamma_1)\\}} \\
&\qquad\qquad\qquad
A_{+\infty}(\gamma_1)(y)\prod_{i=1}^n\prod_{j=1}^{l_i}
(-1)^{e_{i,j} -r_{i,j}(g-1)}[e_{i,j} -r_{i,j}(g-1)]_y H(\eta_{i,j})(y)\\
\eal
\ee
where $\gamma=(r_1,e_1)\in \IZ_{\geq 1}\times \IZ$,
$\eta_{i,j}=(e_{i,j}, r_{i,j})\in \IZ_{\geq 1}\times \IZ$, $1\leq i\leq n$, $1\leq j\leq l_i$. T
Moreover, the sum
in the right hand side of equation \eqref{eq:wallBB} is  finite for any fixed
$\gamma=(r,e)$.

Then in equation \eqref{eq:wallA} $A_{\delta_-}(\gamma)=A_{\delta_+}(\wgamma)$
and
\be\label{eq:wallBC}
\bal
& \mathop{\sum_{l\geq 2}} {1\over (l-1)!} \sum_{\substack{\gamma_1+\cdots +\gamma_l =\gamma\\ \mu(\gamma_i)=\mu(\gamma), 1\leq i\leq l}}
A_{\delta_-}(\gamma_1)(y)\prod_{i=2}^{l-1} e^{e_i-r_i(g-1)} [e_i-r_i(g-1)]_y H(\gamma_i)(y)=\\
& \mathop{\sum_{l\geq 2}} {1\over (l-1)!} \sum_{\substack{\gamma_1+\cdots +\gamma_l =\gamma\\ \mu(\gamma_i)=\mu(\gamma), 1\leq i\leq l}}
A_{\delta_+}(\wgamma_1)(y)\prod_{i=2}^{l-1} e^{e_i-r_i(g-1)} [e_i-r_i(g-1)]_y H(\gamma_i)(y)=\\
& \mathop{\sum_{l\geq 2}} {(-1)^{l-1}\over (l-1)!}
\sum_{\substack{\gamma_1+\cdots +\gamma_l =\wgamma\\ \mu(\gamma_i)=\mu(\wgamma), 1\leq i\leq l}}
A_{\delta_+}(\gamma_1)(y)\prod_{i=2}^{l-1} e^{e_i-r_i(g-1)} [e_i-r_i(g-1)]_y H(\gamma_i)(y)=\\
\eal
\ee
by a redefinition of variables.
Substituting \eqref{eq:wallBB} and \eqref{eq:wallBC} in equation \eqref{eq:wallA},
 equation \eqref{eq:higgsrecursionA} follows by simple combinatorics.

\section{Asymptotic refined ADHM invariants from gauge theory}\label{gaugesect}
The main goal of this section is to present a string theoretic derivation
of conjecture (\ref{refasympinvI}).  Readers who are not interested
in this derivation are encouraged to skip this section.

Conjecture (\ref{refasympinvI}) will be shown to follow from type IIA/M-theory duality
using arguments analogous to \cite{Lawrence:1997jr,EK-I,Nekrasov:2002qd,IKP-I,IKP-II,EK-II,Hollowood:2003cv,Konishi-I,LLZ,IKV}.
Summarizing these results, the topological string amplitudes of certain
toric Calabi-Yau threefolds (as well as some nontoric configurations of local rational curves)
were identified with the instanton partition function of five
dimensional gauge theory
compactified on a circle of finite radius. The later has been identified in
\cite{Nekrasov:2002qd}  with the generating function for the
equivariant Hirzebruch genus of the moduli space of torsion free
framed sheaves on the projective plane. A mathematical exposition
can be found for example in  \cite{instcountA,instcountB}.
The relation between topological
strings and five dimensional gauge theory
has been subsequently refined in \cite{IKV}. Moreover, the
refined topological
string partition function constructed in \cite{IKV} has been
conjecturally identified in \cite{DG} with the generation function of refined
Donaldson-Thomas invariants. The present problem requires a version
of this identification for higher genus local curves.

\subsection{Geometric engineering via local ruled surfaces}\label{geomeng}
Working under the same assumptions as in section (\ref{hitchinpairs}),
$M_1, M_2$ are line bundles on the curve $X$ so that $M_1\otimes_X M_2\simeq K_X^{-1}$, $p=d(M_1)\geq 0$ and $M_1\simeq \CO_X$ if $p=0$.
Let $Y$ be the total space of the rank two vector bundle
$M_1^{-1}\oplus M_2^{-1}$ on $X$, which is a noncompact
Calabi-Yau threefold under the current assumptions.
There is a torus action ${\bf S}\times Y \to Y$ scaling $M_1^{-1}$, $M^{-1}_2$
with characters $t$, $t^{-1}$, so that $Y$ is equivariantly K-trivial.
In principle the relevant five dimensional gauge theory should be constructed
by geometric engineering, that is identifying the low energy effective action
of an M-theory supersymmetric background defined by $S^1\times Y$. This direct approach is somewhat problematic in the present case. A much clearer picture emerges considering a
different local Calabi-Yau threefold constructed as follows.

Let $S$ be the total space of the projective bundle
$\IP(\CO_X\oplus M_1)$. $S$ is a smooth geometrically ruled surface over
$X$ and it has two canonical sections $X_1$, $X_2$ with normal
bundles
\[
N_{X_1/S} \simeq M_1^{-1}, \qquad N_{X_2/S}\simeq M_{1}
\]
respectively. Note that the cone of effective curve classes on $S$ is generated
by the section class $[X_2]$ and the fiber class.

Let $Z$ be the total space of the canonical bundle $K_S$, which is again a noncompact
Calabi-Yau threefold. The normal bundle to $X_1$ in $Z$
is
\[
N_{X_1/Z} \simeq M_1^{-1} \oplus K_X\otimes_X M_1 \simeq M_1^{-1}\otimes M_2^{-1},
\]
therefore the total space of $N_{X_1/Z}$ is isomorphic to $Y$.
Moreover, there is a torus action ${\bf S}\times Z\to Z$
so that $Z$ is equivariantly Calabi-Yau and the induced
torus action on $N_{X_1/Z}$ is compatible with the torus action on $Y$.

Now the main observation is that the local threefold $Z$ engineers a supersymmetric
five dimensional $SU(2)$ gauge theory with $g$ adjoint hypermultiplets on
$\IC^2\times S^1$, where $g$ is
the genus of $X$ \cite{Katz:1996ht}. The integer $p=\mathrm{deg}(M_1)$ corresponds to the
level of the five dimensional Chern-Simons term \cite{Tachikawa:2004ur}.
Therefore by analogy with \cite{Lawrence:1997jr,EK-I,Nekrasov:2002qd,IKP-I,IKP-II,EK-II,Hollowood:2003cv,Konishi-I,LLZ,IKV},
the refined topological string partition function  of $Z$ should be related with the
equivariant instanton partition function ${\CZ}_{inst}^{(p)}(Q,\epsilon_1,\epsilon_2,
a_1,a_2,y)$,
which has been constructed
in \cite{Nekrasov:2002qd}.
As explained in detail in the next subsection,  ${\CZ}_{inst}^{(p)}(Q,\epsilon_1,\epsilon_2,
a_1,a_2,y)$ is the generating function for the $\chi_y$-genus of a certain
holomorphic bundle on a partial compactification of the instanton moduli space.
In particular $\epsilon_1,\epsilon_2$, $a_1,a_2$ are equivariant parameters for a
natural
torus action, $Q$ is a formal
variable counting instanton charge, and $y$ is another formal variable.

In order to make string duality predictions more precise, let $Q_f,Q_b$ be formal
symbols associated to the fiber class, respectively section class $[X_1]$ on $Z$.
Then string duality predicts that there is a factorization
\be\label{eq:partfctfact}
{\CZ}_{ref}(Z; Q_f,Q_b,q,y) = {\CZ}_{ref}^{pert}(Z;Q_f,q,y)
{\CZ}_{ref}^{nonpert}(Z; Q_f,Q_b,q,y)
\ee
into a perturbative, respectively nonperturbative parts. Moreover,
and there is an  identification
\[
{\CZ}_{ref}^{nonpert}(Z;
Q_f,Q_b,q,y) = {\CZ}_{inst}^{(p)}(Q,\epsilon_1,\epsilon_2,a_1,a_2 ,y)
\]
subject to certain
duality relations between the formal parameters in the two partition functions.

Next note that only non-negative powers of $Q_b, Q_f$
can appear in  ${\CZ}_{ref}(Z; Q_f,Q_b,q,y)$ since
the section class $[X_1]$ and the fiber class generate the Mori cone of
$S$. Similarly, only non-negative powers of $Q_f$ can appear in
${\CZ}_{ref}^{pert}(Z;Q_f,q,y)$, which represents the contribution of
pure fiber classes to ${\CZ}_{ref}(Z; Q_f,Q_b,q,y)$.
Therefore ${\CZ}_{ref}(Z; Q_f,Q_b,q,y)$, ${\CZ}_{ref}^{pert}(Z;Q_f,q,y)$
have
well defined  specialization at $Q_f=0$.
Moreover, by construction ${\CZ}_{ref}^{pert}(Z;Q_f,q,y)\big|_{Q_f=0}=1$.
Therefore ${\CZ}_{ref}^{nonpert}(Z; Q_f,Q_b,q,y)$
has well defined specialization at $Q_f=0$ as well,
which is determined by the instanton expansion
$\CZ^{(p)}_{inst}(Q,\epsilon_1,\epsilon_2,a_1,a_2,y)$.
The refined theory of the local threefold $Y$ is then determined by identifying the
contributions of curves supported on the section $X_1$ to
${\CZ}_{ref}^{nonpert}(Z; Q_f,Q_b,q,y)\big|_{Q_f=0}$.
 Computations will be carried out in detail in the next subsections,
resulting in explicit formulas for the instanton partition function and duality
relations among formal variables.

\subsection{Hirzebruch genus}
Let $M(r,k)$ denote the moduli space of rank $r$
framed torsion-free sheaves $(F,f)$ on $\IP^2$
with second Chern class $k\in \IZ_{\geq 0}$.
The framing data is an isomorphism
\be\label{eq:framingA}
f : F|_{\IP^1_\infty} \to \CO_{\IP^1_\infty}^{\oplus r}.
\ee
$M(r,k)$ is a smooth quasi-projective fine moduli space i.e. there is an
universal framed sheaf $({\sf F}, {\sf f})$ on $M(r,k)\times \IP^2$.
Let ${\sf V}= R^1p_{1*}{\sf F}\otimes p_2^*\CO_{\IP^2}(-1)$ where
$p_1,p_2: M(r,k)\times \IP^2\to M(r,k), \IP^2$ denote the canonical
projections. It follows from \cite{hilblect} that ${\sf V}$ is a locally free
sheaf of rank $k$ on $M(r,k)$.

There is a torus ${\bf T} = \IC^\times \times \IC^\times\times (\IC^\times)^{\times r}$ action on
acting on $M(r,k)$, where the action of the first two factors is induced by the canonical
action on $\IC^\times\times \IC^\times $ on $\IP^2$, and the last $r$ factors
act linearly on the framing.
According to \cite{instcountA} the
 fixed points of the {\bf T}-action on $M(r,k)$ are isolated and classified by
collections of Young diagrams ${\underline Y}=(Y_1,\ldots, Y_r)$
so that the total number of
boxes in all diagrams is $|{\underline Y}|=|Y_1|+\cdots |Y_r| =k$.
Let $\CY_{r,k}$ denote the set of all such $r$-uples of Young diagrams.
Note also that both the holomorphic cotangent bundle $T_{M(r,k)}^\vee$ and
the bundle ${\sf V}$ constructed in the previous paragraph carry canonical
equivariant structures.

The K-theoretic instanton partition function of an $SU(2)$ theory with $g$ adjoint
hypermultiplets and a level $p$ Chern-Simons term is given by the equivariant residual Hirzebruch genus of the holomorphic
{\bf T}-equivariant bundle
\[
(T^\vee_{M(2,k)})^{\oplus g} \otimes (\mathrm{det}\, {\sf V})^{-p}.
\]
This is defined by  equivariant localization as follows \cite{instcountB,LLZ}.
Let $(\epsilon_1,\epsilon_2,a_1,a_2)$ be equivariant parameters associated to the torus {\bf T}.
Then the localization formula yields \cite{instcountB,LLZ}
\be\label{eq:instpartfctA}
{\CZ}^{(g,p)}_{inst}(Q,\epsilon_1,\epsilon_2,a_1,a_2,y) = \mathop{\sum_{k=0}^\infty}_{}
Q^k{\CZ}^{(g,p)}_{k}(\epsilon_1,\epsilon_2,a_1,a_2; y)
\ee
where ${\CZ}^{(g,p)}_{0}(\epsilon_1,\epsilon_2,a_1,a_2; y) =1$ and
\be\label{eq:instpartfctB}
\bal
{\CZ}^{(g,p)}_{k}(\epsilon_1,\epsilon_2,a_1,a_2; y)
=
\sum_{\uY\in \CY_{2,k}}\
\prod_{\alpha=1}^2 & \bigg(e^{-|Y_\alpha| a_\alpha}\prod_{(i,j)\in Y_\alpha}
e^{(i-1)\epsilon_1+(j-1)\epsilon_2}\bigg)^p\\
\prod_{\alpha,\beta=1}^2& \prod_{(i,j)\in Y_\alpha}
{\left(1- ye^{(Y^t_{\beta,j}-i)\epsilon_1 - (Y_{\alpha,i}-j+1)\epsilon_2 + 
a_{\alpha\beta}}
\right)^g
\over
\left(1- e^{(Y^t_{\beta,j}-i)\epsilon_1 - (Y_{\alpha,i}-j+1)\epsilon_2 + 
a_{\alpha\beta}}
\right)}
\\
& \prod_{(i,j)\in Y_\beta} {\left(1-ye^{-(Y^t_{\alpha,j}-i+1)\epsilon_1+(Y_{\beta_,i}-j)\epsilon_2 + a_{\alpha\beta}}\right)^g \over
\left(1-e^{-(Y^t_{\alpha,j}-i+1)\epsilon_1+(Y_{\beta_,i}-j)\epsilon_2 + 
a_{\alpha\beta}}\right)}
\eal
\ee
where for any Young tableau $Y$, $Y_i$, $i\in \IZ_{\geq 1}$
denotes the length of the $i$-th
column and $Y^t$ denotes the transpose of $Y$.
If $i$ is greater than the number of columns of $Y$, $Y_i=0$.
Moreover $a_{\alpha\beta}=a_\alpha-a_\beta$ for any $\alpha,\beta=1,2$. 

\subsection{Comparison with the ruled vertex}
A conjectural formula for the unrefined topological string partition function
$\CZ_{top}(Z; Q_f,Q_b,q)$ of the threefold $Z$ has been derived from large
$N$ duality in \cite{ruled}. The purpose of this subsection, is to show that
$\CZ_{top}(Z; Q_f,Q_b,q)$ has a factorization of the form \eqref{eq:partfctfact}
and there is an identification
\[
\CZ_{top}^{nonpert}(Z; Q_f,Q_b,q)= {\CZ}^{(g,p)}_{inst}(Q,\epsilon_1,\epsilon_2,a_1,a_2,y)
\]
subject to certain duality relations between the formal parameters.
This will be a confirmation of duality predictions for local ruled surfaces
in the unrefined case.
Moreover, it will provide a starting point for understanding this correspondence
in the refined case.

By analogy with \cite{IKP-I,LLZ}, first set
\be\label{eq:dualityrelA}
-\epsilon_1=\epsilon_2=\hbar, \qquad y=1.
\ee
Then a straightforward computation yields
\be\label{eq:instpartfctE}
\bal
 & \CZ_{2,k}^{(g,p)}(-\hbar,\hbar,a_1,a_2, 1)
= \\
&  \mathop{\sum_{Y_1,Y_2}}_{|Y_1|+|Y_2|=k} e^{-p(|Y_1|a_1+|Y_2|a_2)} \prod_{\alpha=1}^2 \prod_{(i,j)\in Y_\alpha}
e^{p(j-i)\hbar}\left(2\, \mathrm{sinh}\, {\hbar\over 2}(Y_{\alpha,i}+Y^t_{\alpha,j}-i-j+1)\right)^{2(g-1)}
\\
& \qquad\qquad\qquad\qquad\qquad \qquad \ \ \prod_{(i,j)\in Y_1}
\left(2\, \mathrm{sinh}\, {1\over 2}(a_{1,2}+(Y^t_{2,j} + Y_{1,i}-i-j+1)\hbar)
\right)^{2(g-1)}\\
& \qquad\qquad\qquad\qquad\qquad \qquad \ \ \prod_{(i,j)\in Y_2}
\left(2\, \mathrm{sinh}\, {1\over 2}(a_{1,2}-(Y^t_{1,j}+Y_{2,i}-i-j+1)\hbar)
\right)^{2(g-1)}\\
\eal
\ee
Using identity \cite[Lemm. 4.4]{LLZ}, which was conjectured in \cite{IKP-I} and proven in
\cite{EK-II}, it follows that
\be\label{eq:instpartfctEB}
\bal
 & \CZ_{2,k}^{(g,p)}(-\hbar,\hbar,a_1,a_2,1)
= \\
& \mathop{\sum_{Y_1,Y_2}}_{|Y_1|+|Y_2|=k} 2^{8(g-1)(|Y_1|+|Y_2|)}
e^{-p(|Y_1|a_1+|Y_2|a_2)} e^{p(\kappa(Y_1)+\kappa(Y_2))\hbar/2} \\
&\qquad \qquad
\prod_{\alpha,\beta=1}^2\mathop{\prod_{i,j=1}^\infty}_{}
\left({\mathrm{sinh}\, {1\over 2}(a_{\alpha,\beta} +(Y_{\alpha,i}-Y_{\beta,j}+j-i)\hbar)\over
 \mathrm{sinh}\, {1 \over 2}(a_{\alpha,\beta}+(j-i)\hbar)} \right)^{2(1-g)}\\
\eal
\ee
where for any Young diagram $Y$
\[
\kappa(Y) = 2\mathop{\sum_{(i,j)\in Y}}(j-i) = |Y| +  \mathop{\sum_{i=1}^{l(Y)}}_{}
(Y_i^2-2iY_i),
\]
$l(Y)$ being the number of rows of $Y$.
Note that $\kappa(Y) = -\kappa(Y^t)$.

The topological
string partition function on $Z$ computed by the ruled vertex formalism
\cite{ruled} is
\be\label{eq:ruledtopstringA}
\bal
\CZ_{top}(Z;q,Q_f,Q_b) & = \sum_{Y_1,Y_2} (K_{Y_1,Y_2}(q,Q_f))^{2(1-g)}
Q_b^{|Y_1|+|Y_2|} Q_f^{p|Y_2|}(-1)^{p(|Y_1|+|Y_2|)}
q^{p(\kappa(Y_2)-\kappa(Y_1))/2} \\
\eal
\ee
where
\[
K_{Y_1,Y_2}(q,Q_f) = \sum_{Y} Q_f^{|Y|}W_{Y_2 Y}(q)W_{Y Y_1}(q)
\]
and
\[
W_{R_1, R_2}(q) = s_{R_2}(q^{-i+1/2})) s_{R_1}(q^{R_{2,i}-i+1/2})
\]
for any two Young tableaus $R_1,R_2$. Here $s_R(x^i)$ denotes the Schur function 
associated to the Young tableau $R$. 

According to \cite{IKP-I,EK-II}, \cite[Thm. 7.1]{LLZ}, $K_{Y_1,Y_2}(q,Q_f)= 
K_{Y_2,Y_1}(q,Q_f)$ and 
\be\label{eq:Kidentity}
{K_{Y_1,Y_2^t}(e^{-z}, e^{-b})\over K_{\emptyset,\emptyset}(e^{-z}, e^{-b})}=
(2^{-4}Q_f^{-1/2})^{|Y_1|+|Y_2|} \prod_{\alpha, \beta=1}^2
\prod_{i,j=1}^\infty {\mathrm{sinh}\ {1\over 2}(b_{\alpha,\beta} +
(Y_{\alpha,i}-Y_{\beta,j}+j-i)z) \over
\mathrm{sinh}\, {1\over 2}(b_{\alpha,\beta} + (j-i)z)}
\ee
where $b_{1,2}=-b_{2,1}=b$. 
Therefore \eqref{eq:ruledtopstringA} is equivalent to 
\be\label{eq:ruledtopstringB}
\bal
\CZ_{top}(Z;q,Q_f,Q_b) 
& = \sum_{Y_1,Y_2} (K_{Y_2 ,Y_1^t}(q,Q_f))^{2(1-g)}
Q_b^{|Y_1|+|Y_2|} Q_f^{p|Y_2|}(-1)^{p(|Y_1|+|Y_2|)} q^{p(\kappa(Y_1)+\kappa(Y_2))/2} \\
\eal
\ee
Setting 
\[
\CZ^{pert}_{top}(Z;q,Q_f,Q_b) = K_{\emptyset,\emptyset}(q,Q_f)^{2(1-g)},\qquad
\CZ^{nonpert}_{top}(Z;q,Q_f,Q_b) = {\CZ_{top}(q,Q_f,Q_b)\over K_{\emptyset,\emptyset}(q,Q_f)^{2(1-g)}}.
\]
identity \eqref{eq:Kidentity} yields 
\be\label{eq:partfctidA}
\CZ^{nonpert}_{top}(Z;q,Q_f,Q_b) = \sum_{k=0}^\infty Q^k 
\CZ^{(g,p)}_{2,k}(-\hbar,\hbar, a_1,a_2; 1)
\ee
for the following change of variables 
\be\label{eq:varchangeA}
Q_f=e^{a_{12}},\qquad q=e^{\hbar},\qquad 
Q=Q_bQ_f^{g-1}, \qquad e^{a_1}=-1.
\ee
This is a concrete confirmation of duality predictions in the unrefined case.
The refined case is the subject of the next subsection.

\subsection{Refinement}

As explained at the end of subsection (\ref{geomeng}), string duality predicts
that the nonperturbative part of the refined topological partition function of $Z$
is determined by instanton partition function $\CZ_{inst}^{(p)}(Q,\epsilon_1,
\epsilon_2,a_1, a_2,y)$ provided one finds the correct identification of formal
parameters as in \cite{Hollowood:2003cv,IKV}.
Although local ruled surfaces are not discussed in \cite{Hollowood:2003cv,IKV},
a careful inspection of the cases discussed there leads to the following construction.

Recall that the   contribution
of a fixed point $(Y_1,Y_2)\in \CY_{2,k}$ for some arbitrary $k\geq 1$
to the right hand side of the localization
formula \eqref{eq:instpartfctB} is
\be\label{eq:locfixedpt}
\bal
\prod_{\alpha=1}^2 & \bigg(e^{-|Y_\alpha| a_\alpha}\prod_{(i,j)\in Y_\alpha}
e^{(i-1)\epsilon_1+(j-1)\epsilon_2}\bigg)^p\\
\prod_{\alpha,\beta=1}^2& \prod_{(i,j)\in Y_\alpha}
{\left(1- ye^{(Y^t_{\beta,j}-i)\epsilon_1 - (Y_{\alpha,i}-j+1)\epsilon_2 + a_\alpha-
a_\beta}
\right)^g
\over
\left(1- e^{(Y^t_{\beta,j}-i)\epsilon_1 - (Y_{\alpha,i}-j+1)\epsilon_2 + a_\alpha-
a_\beta}
\right)}
\\
& \prod_{(i,j)\in Y_\beta} {\left(1-ye^{-(Y^t_{\alpha,j}-i+1)\epsilon_1+(Y_{\beta_,i}-j)\epsilon_2 + a_\alpha-a_\beta}\right)^g \over
\left(1-e^{-(Y^t_{\alpha,j}-i+1)\epsilon_1+(Y_{\beta_,i}-j)\epsilon_2 + a_\alpha-
a_\beta}\right)}
\eal
\ee
Let $\CZ^{(g,p)}_{(\emptyset,Y)}(q_1,q_2,Q_f,y)$ be the expression obtained by setting
$q_1=e^{-\epsilon_1}$, $q_2=e^{-\epsilon_2}$ and
\[
Q_f= e^{a_{12}},\qquad e^{a_1}=-1
\]
in \eqref{eq:locfixedpt}.
Note that a simple power counting argument shows that the expression
\[
Q_f^{(g-1)|Y|}
\CZ^{(g,p)}_{(Y,\emptyset)}(q_1,q_2,Q_f,y)
\]
has well defined specialization $\CZ^{(g,p)}_{(Y,\emptyset)}(q_1,q_2,y)^{(0)}$
at $Q_f=0$, for any $Y$. 
Then, for any $r\in \IZ_{\geq 1}$, any Young diagram $Y$ with $|Y|=r$,
and any $p\in \IZ$ let
\be\label{eq:refD}
\bal
\Omega^{(g,p)}_{Y}(\lambda,y) = y^{2|Y|}\lambda^{(g-1)|Y|}
\CZ^{(g,p)}_{(Y,\emptyset)}(\lambda^{-1} y,\lambda y,y^{-1})^{(0)}.
\eal
\ee
Then string duality predicts that the generating function of asymptotic
singly refined ADHM invariants
is given by
\[
\CZ_{+\infty}(\CX,r;\lambda,y) = \sum_{\substack{|Y|=r}} \Omega^{(g,p)}_Y(\lambda,y).
\]
Formula \eqref{eq:omegaformulaI} follows by a straightforward computation.

\subsection{Double Refinement} 
Physical arguments \cite{DM-crossing} present compelling evidence for the existence 
of a doubly refined BPS counting function, which is graded by $U(1)_R$ charge in addition 
to spin quantum number. In this section it is conjectured that the doubly refined partition function 
of asymptotic ADHM invariants is obtained again from the equivariant instanton sum 
\eqref{eq:instpartfctB} by a different specialization of the equivariant parameters. 
Namely, for $r\in \IZ_{\geq 1}$, any Young diagram $Y$ with $|Y|=r$,
and any $p\in \IZ$ let
\be\label{eq:refX} 
\Omega^{(g,p)}_{(Y,\emptyset)}(\lambda, u^{1/2},v^{1/2}) = 
u^{(g+1)|\mu|}v^{(g-1)|\mu|}
\CZ_{(Y,\emptyset)}^{(g,p)}(\lambda^{-1} (uv)^{1/2},
 \lambda (uv)^{1/2},u^{-1})^{(0)}.
\ee
The generating function of doubly refined asymptotic ADHM invariants is then conjectured 
to be 
\[
\CZ_{+\infty}(\CX,r;\lambda,u,v) = \sum_{\substack{|Y|=r}} \Omega^{(g,p)}_Y(\lambda, u^{1/2},v^{1/2}).
\]
A straightforward computation yields formula \eqref{eq:omegaformulaII}.
In conjunction with the doubly refined wallcrossing conjecture (\ref{bigradedinv}), 
the above formula will be shown to yield correct results for the Hodge polynomial of the Hitchin moduli
space in many examples recorded in appendix (\ref{hodgeapp}).

\subsection{Localization interpretation for $r=2$}
Suppose the conditions of section (\ref{hitchinpairs}) are satisfied, that is
$p\geq 0$, and  $M_1=\CO_X$, $M_2=K_X^{-1}$ if $p=0$. The goal of this
section is to discuss the geometric interpretation of conjecture (\ref{refasympinvI})
for $r=1,2$. The main observation is that in these cases, equation \eqref{eq:refC}
can be interpreted as a sum of contributions of torus fixed loci in the moduli
space $\mfm_{+\infty}^{ss}(\CX,r,e)$. However, a rigorous geometric computation
would require a localization theorem for the refined Donaldson-Thomas invariants
defined in \cite{wallcrossing}, which has not been been formulated and proven so far.

First let $r=1$. The moduli stack of $\delta$-semistable ADHM sheaves of type $(1,e)$
on $X$
with $\delta>0$ and $e\geq 0$ is a $\IC^\times$-gerbe over the smooth
variety
\be\label{eq:modspaceA}
 S^e(X) \times  H^0(X, M_1^{-1})\times H^0(X,M_2^{-1}).
\ee
A $\IC$-valued point of $\mfm_{\delta}^{ss}(\CX,1,e)$ is an ADHM sheaf of the form
$(E, \Phi_1,\Phi_2, 0, \psi)$ where $E$ is a degree $e$ line bundle on $X$,
$\Phi_1\in \mathrm{Hom}_X(E\otimes_X M_1,E)\simeq H^0(X, M_1^{-1})$, $\Phi_2\in
\mathrm{Hom}_X(E\otimes_X M_2, E)\simeq H^0(X,M_2^{-1})$ and
$\psi\in H^0(X,E)$. The $\delta$-stability condition, $\delta>0$ is equivalent to
$\psi$ not identically zero.   Obviously, the moduli stack is empty if $e<0$.

The fixed point conditions require $\Phi_1=0$, $\Phi_2=0$. Therefore
the torus fixed locus  is a $\IC^\times$-gerbe over the symmetric
product $S^e(X)$.

Conjecture \eqref{refasympinvI} and equation
\eqref{eq:omegaformulaI}     yield
\be\label{eq:asymprankoneA}
\CZ_{+\infty}(\CX,1;\lambda,y) = (-1)^p y^{1-g}\frac{(1-\lambda)^{2g}}{(1-\lambda y)(1-\lambda y^{-1})}.
\ee
Now recall Macdonald's formula
\be\label{eq:macdonald}
\sum_{n\geq 0}P_z(S^n(X)) x^n = {(1-xz)^{2g}\over (1-x)(1-xz^2)}.
\ee
for the generating function of Poincar\'e polynomials of symmetric products of $X$.
Then equations \eqref{eq:asymprankoneA} and \eqref{eq:macdonald} imply
\be\label{eq:asymprankoneB}
\CZ_{+\infty}(\CX,1;\lambda,y) =(-1)^p\sum_{e\geq 0} \lambda^e y^{1-g-e}P_y(S^e(X))
\ee
for all $e\in \IZ_{\geq 0}$.

Next let $r=2$. Property $(B.2)$ implies that the moduli space $\mfm_{+\infty}^{ss}(\CX,r,e)$ is empty unless $e\geq 2-2g$. Assuming this to be the case, a straightforward
analysis shows that the components of the torus fixed locus are of two types.
The ADHM sheaves corresponding to the $\IC$-valued fixed points are presented as follows.

$(i)$ $E\simeq E_{-1}\oplus E_0$, $\Phi_2=0$, $\mathrm{Im}(\psi)\subseteq E_0$
and
\[
\Phi_1=\left[\begin{array}{cc} 0 & \varphi  \\ 0 & 0 \\  \end{array}\right]
\]
with $\varphi: E_0\otimes_X M_1\to E_{-1}$ a nontrivial morphism of line bundles.
Components of this type are isomorphic to $\IC^\times$-gerbes over the smooth varieties
\[
S^{e_0}(X) \times S^{e_{-1}-e_0-p}(X)
\]
where $0\leq e_0 \leq e_{-1}-p$ and $e_0+e_{-1}=e$.

$(ii)$  $E\simeq E_0\oplus E_1$, $\Phi_1=0$, $\mathrm{Im}(\psi)\subseteq E_0$
and
\[
\Phi_2=\left[\begin{array}{cc} 0 & 0 \\ \varphi & 0 \\  \end{array}\right]
\]
with $\varphi: E_0\otimes_X M_2\to E_1$ a nontrivial morphism of line bundles.
Components of this type are isomorphic to $\IC^\times$-gerbes over the smooth varieties
\[
S^{e_0}(X) \times S^{e_1-e_0+2g-2+p}(X)
\]
where $0\leq e_0 \leq e_1+2g-2+p$ and $e_0+e_1=e$.

Note that in both cases, the moduli stack of asymptotically stable ADHM sheaves is
not smooth along the fixed loci, although the fixed loci are smooth.

Conjecture (\ref{refasympinvI}) and equation \eqref{eq:omegaformulaI} yield
\be\label{eq:asympranktwoA}
\bal
\CZ_{+\infty}(\CX,2;\lambda,y) & = \Omega^{(g,p)}_{\tableau{2}}(\lambda, y) +
\Omega^{(g,p)}_{\tableau{1 1}}(\lambda,y)\\
\Omega^{(g,p)}_{\tableau{1 1}}(\lambda, y)& = (\lambda^{-1} y)^{-p}y^{2-2g} \frac{(1-\lambda^2y^{-1})^{2g}(1-\lambda)^{2g}}{(1-\lambda^2)(1-\lambda^2y^{-2})(1-\lambda y)(1-\lambda y^{-1})}\\
\Omega^{(g,p)}_{\tableau{2}}(\lambda,y)&= (\lambda y)^{-p}y^{4-4g}
  \lambda^{2-2g}\frac{(1-\lambda^2 y)^{2g}(1-\lambda)^{2g}}{(1-\lambda^2) (1-\lambda^2 y^2) (1 - \lambda y) (1 - \lambda
  y^{-1})}\\
\eal
\ee
A straightforward computation using equation \eqref{eq:macdonald}
yields
\be\label{eq:asympranktwoB}
\bal
\Omega^{(g,p)}_{\tableau{1 1}}(\lambda, y)& = \mathop{\sum_{e\geq p}}_{} \lambda^e
\mathop{\sum_{e_0+e_{-1}=e}}_{0\leq e_0 \leq e_{-1}-p}
 y^{2-2g-p-e_0}\ y^{-e_0}P_{y}(S^{e_0}(X))\\
 &\qquad   \qquad\qquad  \qquad \quad     y^{-e_{-1}+e_0+p}
 P_{y}(S^{e_{-1}-e_0-p}(X)) \\
\Omega^{(g,p)}_{\tableau{2}}(\lambda,y)&= \mathop{\sum_{e\geq 2-2g-p}}_{}
\lambda^e \mathop{\sum_{e_0+e_{1}=e}}_{0\leq e_0 \leq e_{1}+2g-2+p}
 y^{e_0-p}\ y^{-e_0}P_{y}(S^{e_0}(X))\\
 & \qquad\qquad\qquad\
 y^{-e_1+e_0-2g+2-p}P_{y}(S^{e_1-e_0+2g-2+p}(X))\\
\eal
\ee
Given the explicit description of the fixed loci, equations \eqref{eq:asymprankoneA},
\eqref{eq:asymprankoneB}, \eqref{eq:asympranktwoA}, \eqref{eq:asympranktwoB}
clearly suggest an equivariant localization theorem for refined ADHM invariants.
 Such a formula would presumably allow a rigorous computation of the
polynomial weights assigned to each component of the fixed locus.

For future reference, let us record the expressions $\Omega_Y^{(p)}(\lambda,y)$
for $|Y|=3$.
\be\label{eq:asymprankthree}
\bal
\Omega^{(g,p)}_{\tableau{1 1 1}}(\lambda,y) & = (-1)^p(\lambda^{3}y^{-3})^py^{3-3g}\frac{(1-\lambda)^{2g}(1-\lambda^2y^{-1})^{2g}(1-\lambda^3y^{-2})^{2g}}{(1-\lambda
y)(1-\lambda y^{-1})(1-\lambda^2y^{-2})(1-\lambda^2)(1-\lambda^3y^{-3})(1-\lambda^3y^{-1})}\\
\Omega^{(g,p)}_{\tableau{2 1}}(\lambda,y)  & = (-1)^py^{2p}y^{5-5g}\lambda^{2-2g}\frac{(1-\lambda)^{4g}(1-\lambda^3)^{2g}}{(1-\lambda
y)^2(1-\lambda y^{-1})^2(1-\lambda^3 y)(1-\lambda^3y^{-1})}\\
\Omega^{(g,p)}_{\tableau{3}}(\lambda,y)  & = (-1)^p(\lambda^{-3}
y^{-3})^py^{9-9g}\lambda^{6-6g}\frac{(1-\lambda)^{2g}(1-\lambda^2y)^{2g}(1-\lambda^3y^2)^{2g}}{(1-\lambda
y)(1-\lambda
y^{-1})(1-\lambda^2y^2)(1-\lambda^2)(1-\lambda^3y^3)(1-\lambda^3y)}\\
\eal
\ee

\section{Examples, comparison with existing results}\label{computations}

This section will present several concrete results for Poincar\'e polynomials
of moduli spaces of Hitchin pairs obtained from the recursion relation
\eqref{eq:higgsrecursionA} and conjecture
\ref{refasympinvI}.  In all cases considered below, these results are identical
to the computations of Hitchin \cite{hitchin-selfd} and Gothen
\cite{Bettinumbers}, as well as the conjecture
of Hausel and Rodriguez-Villegas \cite{HRV,Mirror-Hodge}, which are briefly reviewed
in appendix (\ref{survey}).
In addition, entirely analogous computations  have been done for the Hodge polynomial
of moduli spaces of pairs, employing the doubly refined version of the recursion
formula and conjecture (\ref{refasympinvII}). The results are presented in
appendix (\ref{hodgeapp}). Again, all cases considered there are in agreement with
the results of \cite{hitchin-selfd,Bettinumbers,HRV,Mirror-Hodge}.

In order to simplify the formulas
 set $\wA_{+\infty}(r,e)(y) = (-1)^{rp} A_{+\infty}(r,e)(y)$,
$\wH(r,e)(y)=(-1)^{e-r(g-1-p)}H(r,e)(y)$ for all $(r,e)\in \IZ_{\geq 1}\times \IZ$.
Then equation \eqref{eq:higgsrecursionA} becomes
\be\label{eq:higgsrecursionC}
\bal
& [e-r(g-1)]_y \wH(\gamma)(y) =
\wA_{+\infty}({\gamma})(y)-\wA_{+\infty}(\wgamma)(y)\\
& +
\mathop{\sum_{l\geq 2}}_{}{(-1)^{l-1}\over (l-1)!}
\sum_{\substack{ \gamma_1,\ldots,\gamma_l \in \IZ_{\geq 1}\times \IZ\\
\gamma_1+\cdots+\gamma_l=\gamma\\
\mu_0(r)\leq \mu(\gamma)< \mu(\gamma_i), \ 2\leq i\leq l,\\
 \mu_0(r)\leq \mu(\gamma_1)}}
\wA_{+\infty}(\gamma_1)(y)
\prod_{i=2}^l [e_i-r_i(g-1)]_y \wH(r_i,e_i)(y)\\
& -\mathop{\sum_{l\geq 2}}_{}{(-1)^{l-1}\over (l-1)!}
\sum_{\substack{ \gamma_1,\ldots,\gamma_l \in \IZ_{\geq 1}\times \IZ\\
\gamma_1+\cdots+\gamma_l=\wgamma\\
\mu_0(r)\leq  \mu(\wgamma)\leq \mu(\gamma_i), \ 2\leq i\leq l\\
\mu_0(r)\leq \mu(\gamma_1)}}
\wA_{+\infty}(\gamma_1)(y)
\prod_{i=2}^l[e_i-r_i(g-1)]_y \wH(r_i,e_i)(y)\\
& -\mathop{\sum_{l\geq 2}}_{} {1\over l!}
\sum_{\substack{
\gamma_1,\ldots,\gamma_l \in \IZ_{\geq 1}\times \IZ\\
\gamma_1+\cdots+\gamma_l=\gamma\\
\mu(\gamma)=\mu(\gamma_i),\ 1\leq i\leq l\\}}
\prod_{i=1}^l [e_i-r_i(g-1)]_y \wH(r_i,e_i)(y)\\
\eal
\ee
where $\mu_0(r) = -(r-1)(2g-2+p)$, and the sum in the right hand side of equation
\eqref{eq:higgsrecursionA} is finite.

\subsection{Rank $r=1$}
There are no positive critical parameters of type $(1,e)$ for any $e\in \IZ_{\geq 0}$
 The  wallcrossing formula \eqref{eq:wallA} at $\delta_c=0$ reads
\be\label{eq:rankonewallA}
\wA_{+\infty}(1,e)-\wA_{+\infty}(1,-e+2(g-1)) = [e-g+1]_y \wH(1,e).
\ee
Expanding the right hand side of equation \eqref{eq:asymprankoneB} in powers of
$\lambda$ yields
\[
\bal
\wA_{+\infty}(1,e) = y^{1-g}\mathop{\sum_{0\leq k\leq 2g}}_{m,l\geq 0,\
k+l+m=e} (2g,k) (-1)^k y^{l-m}
\eal
\]
for any $e\geq 0$, where $(2g,k) = {(2g)!\over k! (2g-k)!}$ are binomial
coefficients.
A series of elementary manipulations further yield
\[
\bal
\wA_{+\infty}(1,e) & = y^{1-g} \mathop{\sum_{0\leq k\leq 2g}}_{l\geq 0, \
l+k \leq e} (2g,k) (-1)^k y^{2l+k-e}\\
& = y^{1-g} \mathop{\sum_{0\leq k\leq 2g}}_{l\geq 0, \
l+k \leq e} (2g,k) (-1)^k y^{k-e} {1-y^{2e-2k+2}\over 1-y^2}\\
& = {y^{1-g}\over 1-y^2}\mathop{\sum_{0\leq k\leq 2g}}_{
k \leq e} (2g,k) (-1)^k  \big(y^{k-e}-y^{e-k+2}\big)\\
\eal
\]
for any $e\geq 0$. In order to compute the left hand side of equation
\eqref{eq:rankonewallA}, it is convenient to consider three cases.

$a)\ 0\leq e\leq 2g-2$. Then
\[
\bal
& \wA_{+\infty}(1,e) -\wA_{+\infty}(1,-e+2(g-1))   = \\
& {y^{1-g}\over 1- y^2} \bigg[\mathop{\sum_{k=0}^e}_{} (2g,k)
(-1)^k y^{k-e}+ \mathop{\sum_{k=0}^{2g-2-e}}_{} (2g,k)
(-1)^k y^{2g-2-e-k}\bigg] \\
& -{y^{1-g}\over 1- y^2} \bigg[\mathop{\sum_{k=0}^e}_{} (2g,k)
(-1)^k y^{e-k+2}+ \mathop{\sum_{k=0}^{2g-2-e}}_{} (2g,k)
(-1)^k y^{k+e-2g+2}\bigg] \\
& = {y^{1-g}\over 1- y^2}y^{-e} \bigg[\mathop{\sum_{k=0}^e}_{} (2g,k)
(-1)^k y^{k}+ \mathop{\sum_{k=e+2}^{2g-2}}_{} (2g,k)
(-1)^k y^{k}\bigg] \\
&\ \ \ -{y^{1-g}\over 1- y^2}y^{e+2}
\bigg[\mathop{\sum_{k=0}^e}_{} (2g,k)
(-1)^k y^{-k}+ \mathop{\sum_{k=e+2}^{2g-2}}_{} (2g,k)
(-1)^k y^{-k}\bigg] \\
& = -{y^{1-g}\over 1- y^2} \bigg[y^{e+2}(1-y^{-1})^{2g}-y^{-e}(1-y)^{2g}\bigg]
=  {y^{e-g+1}-y^{-e+g-1}\over y - y^{-1}} {(1-y)^{2g}\over y^{2g-1}}\\
\eal
\]

$b)\ e=2g-1$. Then $\wA_{+\infty}(-e+2g-2)=0$ and
\[
\bal
\wA_{+\infty}(1,2g-1) & = {y^{1-g}\over 1-y^2}
\mathop{\sum_{k=0}^{2g-1}} (2g,k) (-1)^k \big( y^{k-2g+1}-y^{2g-k+1}\big) \\
& = {y^{1-g}\over 1-y^2} \bigg[y^{1-2g}(1-y)^{2g} - y^{2g+1}(1-y^{-1})^{2g}\bigg]
= {y^{g}-y^{-g}\over y - y^{-1}} {(1-y)^{2g}\over y^{2g-1}}\\
\eal
\]

$c)\ e\geq 2g$. Then $\wA_{+\infty}(-e+2g-2)=0$ and a similar computation yields
\[
\wA_{+\infty}(1,e) = {y^{e-g+1}-y^{-e+g-1}\over y - y^{-1}} {(1-y)^{2g}\over y^{2g-1}}.
\]
In conclusion,
\be\label{eq:rankoneB}
\wH(1,e)(y) = {(1-y)^{2g}\over y^{2g-1}}
\ee
for all $e\geq 0$, hence also for all $e\in \IZ$.

The moduli space of rank one semistable Hitchin pairs of any degree $e\in \IZ$
is isomorphic to
\[
H^0(X, M_1^{-1})\times H^0(X, M_2^{-1}) \times J_e(X)
\]
where $J_e(X)$ is the degree $e$ Jacobian of $X$. Obviously the formula
\eqref{eq:rankoneB} can be rewritten as
\[
  \wH(1,e)(y) = y^{1-2g}P_y(J_e(X))
    \]
  for any $e\in \IZ$.

\subsection{Rank $r=2$}
According to property (B.2) in section (\ref{review}), all invariants $A_\delta(1,e)(y)$ are
zero for $e<0$. It will be convenient to distinguish two cases, depending on the parity of
$e$. By convention, any sum in the following formulas is zero if the lower summation bound exceeds the upper summation
bound.

$a)\ e=2n,\ n\in \IZ$.
Then equation \eqref{eq:higgsrecursionC} reduces to
\be\label{eq:ranktworecA}
\bal
& [2n-2g+2]_y \wH(2,2n)(y) =\wA_{+\infty}(2,2n)- \wA_{+\infty}(2,-2n+4g-4)\\
& -\sum_{e_1=0}^{n-1} \wA_{+\infty}(1,e_1)(y)[2n-e_1-g+1]_y \wH(1,2n-e_1)(y)\\
& +\sum_{e_1=0}^{2g-2-n}
\wA_{+\infty}(1,e_1)(y) [3g-3-2n-e_1]_y \wH(1,4g-4-2n)(y) \\
&
-{1\over 2} [n-g+1]_y^2 \wH(1,n)(y)^2.\\
\eal
\ee

$b)\ e=2n+1,\ n\in \IZ$.
Then equation \eqref{eq:higgsrecursionA} reduces to
\be\label{eq:ranktworecB}
\bal
& [2n-2g+3]_y \wH(2,2n+1)(y) = \wA_{+\infty}(2,2n+1) -\wA_{+\infty}(2,4g-5-2n) \\
& -\sum_{e_1=0}^n \wA_{+\infty}(1,e_1)(y)[2n-e_1-g+2]_y \wH(1,2n+1-e_1)(y)\\
& +\sum_{e_1=0}^{2g-3-n}
\wA_{+\infty}(1,e_1)(y) [3g-4-2n-e_1]_y H(1,4g-4-2n-e_1)(y) \\
\eal
\ee
Some concrete results are recorded below. $\wH^{(p)}(r,e)$ denotes the refined Higgs invariant of type $(r,e)$
with coefficient bundles $(M_1,M_2)$ of degrees $(p,2-2g-p)$, $p\geq 0$.
Under the current assumptions, $M_1\simeq \CO_X$ if $p=0$.

\bigskip
\framebox[1.2cm][c]{$ {g=2}$} \par
\bigskip

\[
\wH^{(0)}(2,1)(y)=\frac{(1-y)^4(1+y^2)(1-4y^3+2y^4)}{y^9}
\]
\[
\wH^{(0)}(2,0)(y)=\frac{(1-y)^4(2 + 4 y^2 - 8 y^3 + 7
y^4 - 12 y^5 + 14 y^6 - 4 y^7 + 5 y^8)}{2y^9(1+y^2)}
\]

\[  \wH^{(1)}(2,1)(y)=\frac{(1-y)^4 \left(2 y^8-4 y^7+8 y^6-4 y^5+2 y^4-4
y^3+y^2+1\right)}{y^{11}}\]

\[ \wH^{(1)}(2,0)(y)=\frac{(1-y)^4 \left(3 y^{10}-12 y^9+14 y^8-20 y^7+19 y^6-16 y^5+6 y^4-8 y^3+4 y^2+2\right)}{2 y^{11}
\left(y^2+1\right)}\]

\[ \wH^{(2)}(2,1)(y)=\frac{(1-y)^4 \left(y^2+1\right) \left(2 y^8-8 y^7+6 y^6+2 y^4-4
y^3+1\right)}{y^{13}}\]

\[ \wH^{(2)}(2,0)(y)=\frac{(1-y)^4 \left(5 y^{12}-12 y^{11}+26 y^{10}-28 y^9+33 y^8-24 y^7+20 y^6-16 y^5+6 y^4-8 y^3+4 y^2+2\right)}{2 y^{13}
\left(1+y^2\right)}\]

\bigskip
\framebox[1.2cm][c]{$ {g=3}$} \par
\bigskip

\[
\wH^{(0)}(2,1)(y)= \frac{(1-y)^6}{y^{17}} (1 + y^2 - 6 y^3 + 2 y^4 - 6 y^5 +
17 y^6 - 12 y^7 +18 y^8 - 32 y^9 + 18 y^{10} - 12 y^{11} + 3 y^{12})
\]

\[
\bal
\wH^{(0)}(2,0)(y)=\frac{(1-y)^6}{2y^{17}(1+y^2)} &(2 + 4 y^2 - 12 y^3 + 6 y^4 - 24 y^5 + 38 y^6 - 36 y^7 +
   71 y^8 - 82 y^9 + 87 y^{10} - 68 y^{11} \\
   &+ 57 y^{12} - 18 y^{13} + 7
   y^{14})
   \eal
   \]

   \[ \wH^{(1)}(2,1)(y)=\frac{(1-y)^6 \left(y^2+1\right) \left(3 y^{12}-12 y^{11}+30 y^{10}-20 y^9+3 y^8-12 y^7+15 y^6+2 y^4-6
y^3+1\right)}{y^{19}}\]

\[ \bal \wH^{(1)}(2,0)(y)=\frac{(1-y)^6}{2y^{19}(1+y^2)}&(5 y^{16}-30 y^{15}+57 y^{14}-108 y^{13}+117 y^{12}-134 y^{11}\\&+101 y^{10}-88 y^9+70 y^8-36 y^7+38 y^6-24 y^5+6 y^4-12 y^3+4
y^2+2)\eal\]

     \[ \bal \wH^{(2)}(2,1)(y)=\frac{(1-y)^6}{y^{21}}&(3 y^{16}-18 y^{15}+33 y^{14}-52 y^{13}+48 y^{12}-38 y^{11}\\& +33 y^{10}-32 y^9+18 y^8-12 y^7+17 y^6-6 y^5+2 y^4-6 y^3+y^2+1) \eal \]

   \[ \bal \wH^{(2)}(2,0)(y)=\frac{(1-y)^6}{2y^{21}(1+y^2)}&(7 y^{18}-30 y^{17}+87 y^{16}-120 y^{15}+177 y^{14}-174 y^{13}\\& +163 y^{12}-140 y^{11}+102 y^{10}-88 y^9+70 y^8-36 y^7\\&+38 y^6-24 y^5+6 y^4-12
   y^3+4 y^2+2) \eal\]

   \bigskip
\framebox[1.2cm][c]{$ {g=4}$} \par
\bigskip

      \[
      \bal
      \wH^{(0)}(2,1)(y)=\frac{(1-y)^8}{y^{25}} &(1 + y^2) (1 - 8 y^3 + 2 y^4 + 28 y^6
      - 16 y^7 + 3 y^8 -
   56 y^9 + 56 y^{10} - 24 y^{11}\\ & + 74 y^{12} - 112 y^{13} + 56 y^{14} -
   24 y^{15} + 4 y^{16})
   \eal
   \]

\[
\bal
\wH^{(0)}(2,0)(y)=\frac{(1-y)^8}{2y^{25}(1+y^2)} &(2 + 4 y^2 - 16 y^3 + 6 y^4 - 32 y^5 + 64 y^6 - 48 y^7 + 122 y^8 -
  176 y^9\\& + 180 y^{10} - 304 y^{11} + 379 y^{12} - 424 y^{13} + 548 y^{14} -
  488 y^{15} + 450 y^{16} - 264 y^{17}\\& + 156 y^{18} - 40 y^{19} + 9
  y^{20})
  \eal
  \]

\bigskip
\framebox[1.2cm][c]{$ {g=5}$} \par
\bigskip

 \[
 \bal
 \wH^{(0)}(2,1)(y)=\frac{(1-y)^{10}}{y^{33}} &(1 + y^2 - 10 y^3 + 2 y^4 - 10 y^5 + 47 y^6 - 20 y^7 +
   48 y^8 - 140 y^9 + 93 y^{10} - 150 y^{11}\\& + 304 y^{12} - 270 y^{13} +
   349 y^{14} - 532 y^{15} + 560 y^{16} - 652 y^{17} + 770 y^{18} - 784 y^{19}\\& +
   560 y^{20} - 400 y^{21} + 140 y^{22} - 40 y^{23} + 5
   y^{24})
   \eal
   \]

  \bigskip
  \noindent
In all the above cases, similar computations also show that the invariants $\wH(2,e)$
depend only on the parity of $e\in \IZ$.
Note also that for even $e$ the rank two refined Higgs invariants are rational functions
of $y$ rather than polynomials in $y^{-1},y$. By analogy with the theory of generalized Donaldson-Thomas invariants \cite{genDTI}, this reflects the fact that
in this case the moduli stack ${\mathfrak {Higgs}}^{ss}(\CX,2,e)$ contains
strictly semistable $\IC$-valued points.

\subsection{Rank $r=3$} According to property (B.2) in section (\ref{review}), all invariants $A_\delta(2,e)(y)$ are
zero for $e<2-2g-p$.  Suppose $e=3n+1$,
$n\in \IZ$.
Then equation \eqref{eq:higgsrecursionC} reduces to
\be\label{eq:rankthreerecA}
\bal
& [3n-3g+4]_y \wH(3,3n+1) = \wA_{+\infty}(3,3n+1)(y)-\wA_{+\infty}(3,-3n+6g-7)(y)\\
& -\sum_{e_1=2-2g-p}^{2n} \wA_{+\infty}(2,e_1)[3n+2-g-e_1]_y \wH(1,3n+1-e_1)(y) \\
& - \sum_{e_1=0}^n \wA_{+\infty}(1,e_1) [3n+3-2g-e_1]_y \wH(2,3n+1-e_1)(y)\\
& +{1\over 2} \sum_{e_1=0}^{n-1}\
\sum_{e_2=n+1}^{2n-e_1} \wA_{+\infty}(1,e_1)[e_2-g+1]_y[3n+2-g-e_1-e_2]_y
\wH(1,3n+1-e_1-e_2)(y)^2\\
& +\sum_{e_1=2-2g-p}^{4g-2n-5} \wA_{+\infty}(2,e_1)[5g-6-3n-e_1]_y
\wH(1,6g-7-3n-e_1)(y) \\
& +\sum_{e_1=0}^{2g-n-3} \wA_{+\infty}(1,e_1)
[4g-5-3n-e_1]_y \wH(2,6g-7-3n-e_1)(y)\\
& - {1\over 2}\sum_{e_1=0}^{2g-3-n}\ \sum_{e_2=2g-2-n}^{4g-2n-5-e_1}
\wA_{+\infty}(1,e_1) [e_2-g+1]_y [5g-6-3n-e_1-e_2]_y  \wH(1,6g-7-3n-e_1-e_2)(y)^2 \\
\eal
\ee
Again, some concrete results are recorded below.

\bigskip
\framebox[1.2cm][c]{$ {g=2}$} \par
\bigskip

\[
\bal
\wH^{(0)}(3,1)(y)=&
\frac{(1-y)^4}{y^{19}}\\ &(1 + y^2 - 4 y^3 + 3 y^4 - 8 y^5 + 10 y^6
- 16 y^7 + 29 y^8 - 32 y^9 +
 48 y^{10} - 64 y^{11} \\& + 67 y^{12} - 68 y^{13} + 48 y^{14} - 24 y^{15} + 6
 y^{16})
 \eal
 \]

\[
\bal
\wH^{(1)}(3,1)(y)=\frac{(1-y)^4}{y^{25}}&(6 y^{22}-36 y^{21}+96 y^{20}-168 y^{19}+207 y^{18}-216 y^{17}+210 y^{16}-184 y^{15}\\&+149 y^{14}-120 y^{13}+92 y^{12}-72 y^{11}+49 y^{10}-32
   y^9\\&+29 y^8-16 y^7+10 y^6-8 y^5+3 y^4-4 y^3+y^2+1)
   \eal
   \]

   \[
   \bal
   \wH^{(2)}(3,1)(y)=\frac{(1-y)^4}{y^{31}}&(10 y^{28}-64 y^{27}+184 y^{26}-344 y^{25}+477 y^{24}-560 y^{23}+583 y^{22}\\&-560 y^{21}+522 y^{20}-464 y^{19}+386 y^{18}-320 y^{17}+267
   y^{16}-208 y^{15}\\&+158 y^{14}-124 y^{13}+93 y^{12}-72 y^{11}+49 y^{10}-32 y^9+29 y^8-16 y^7\\&+10 y^6-8 y^5+3 y^4-4
   y^3+y^2+1)
   \eal
   \]

\bigskip
\framebox[1.2cm][c]{$ {g=3}$} \par
\bigskip

\[
\bal
\wH^{(0)}(3,1)(y)=\frac{(1-y)^6}{y^{37}}&(15 y^{32}-120 y^{31}+480 y^{30}-1260 y^{29}+2355 y^{28}-3486 y^{27}\\ &+4189 y^{26}-4416 y^{25}+4315 y^{24}-3922 y^{23}+3399 y^{22}-2860
   y^{21}\\&+2309 y^{20}-1872 y^{19}+1433 y^{18}-1072 y^{17}+861 y^{16}-604 y^{15}\\&+446 y^{14}-336 y^{13}+212 y^{12}-176 y^{11}+105 y^{10}-62 y^9\\&+58 y^8-24 y^7+19
   y^6-12 y^5+3 y^4-6 y^3+y^2+1)
   \eal
   \]

   \[
   \bal
   \wH^{(1)}(3,1)(y)=\frac{(1-y)^6}{y^{43}}&(15 y^{38}-150 y^{37}+690 y^{36}-2010 y^{35}+4110 y^{34}-6542 y^{33}\\&+8598 y^{32}-9930 y^{31}+10427 y^{30}-10254 y^{29}+9672 y^{28}-8800
   y^{27}\\&+7705 y^{26}-6600 y^{25}+5598 y^{24}-4600 y^{23}+3723 y^{22}-3006 y^{21}\\&+2363 y^{20}-1884 y^{19}+1434 y^{18}-1072 y^{17}+861 y^{16}-604 y^{15}+446
   y^{14}\\&-336 y^{13}+212 y^{12}-176 y^{11}+105 y^{10}-62 y^9+58 y^8-24 y^7\\&+19 y^6-12 y^5+3 y^4-6
   y^3+y^2+1)
   \eal
   \]

    \[
    \bal
    \wH^{(2)}(3,1)(y)=\frac{(1-y)^6}{y^{49}}&(21 y^{44}-216 y^{43}+1026 y^{42}-3090 y^{41}+6621 y^{40}-11094 y^{39}\\&+15375 y^{38}-18672 y^{37}+20712 y^{36}-21584 y^{35}+21450 y^{34}-20552
   y^{33}\\&+19178 y^{32}-17460 y^{31}+15503 y^{30}-13546 y^{29}+11706 y^{28}-9952 y^{27}+8316 y^{26}\\&-6912 y^{25}+5736 y^{24}-4650 y^{23}+3741 y^{22}-3012
   y^{21}\\&+2364 y^{20}-1884 y^{19}+1434 y^{18}-1072 y^{17}+861 y^{16}-604 y^{15}\\&+446 y^{14}-336 y^{13}+212 y^{12}-176 y^{11}+105 y^{10}-62 y^9+58 y^8\\&-24 y^7+19
   y^6-12 y^5+3 y^4-6 y^3+y^2+1)
   \eal
   \]

   \bigskip
\framebox[1.2cm][c]{$ {g=4}$} \par
\bigskip

   \[
   \bal
   \wH^{(0)}(3,1)(y)=\frac{(1-y)^8}{y^{55}}&(28 y^{48}-336 y^{47}+2016 y^{46}-7896 y^{45}+22218 y^{44}-48328 y^{43}\\&+84084 y^{42}-122616 y^{41}+155235 y^{40}-176912 y^{39}+186320
   y^{38}-185408 y^{37}\\&+176976 y^{36}-163656 y^{35}+146930 y^{34}-128936 y^{33}+111544 y^{32}-94416 y^{31}\\&+78918 y^{30}-65392 y^{29}+53178 y^{28}-43392
   y^{27}+34620 y^{26}-27288 y^{25}+21936 y^{24}\\&-16728 y^{23}+13005 y^{22}-10064 y^{21}+7290 y^{20}-5760 y^{19}+4077 y^{18}-2880 y^{17}\\&+2278 y^{16}-1416
   y^{15}+1071 y^{14}-744 y^{13}+416 y^{12}-368 y^{11}+185 y^{10}-112 y^9\\&+99 y^8-32 y^7+32 y^6-16 y^5+3 y^4-8
   y^3+y^2+1)
   \eal
   \]

In addition similar computations show that $\wH^{(p)}(3,2)(y)=\wH^{(p)}(3,1)(y)$ in all above examples.

   \appendix
\section{Existing results}\label{survey}

This section is a summary of existing results and conjectures on the cohomology of moduli spaces of Hitchin pairs. The localization computations of Hitchin \cite{hitchin-selfd}
and Gothen \cite{Bettinumbers} as well as the conjectures of
Hausel and Rodriguez-Villegas \cite{HRV,Mirror-Hodge}
will be briefly reviewed. In the first two cases, the localization computations will
be generalized to moduli spaces of Hitchin pairs with coefficient line bundle $L$ of
degree $d(L)=2g-2+p$, $p\geq 0$. As in section (\ref{hitchinpairs}), $L= K_X$
if $p=0$.

For $(r,e)\in \IZ_{\geq 1}\times \IZ$ coprime there is a smooth quasi-projective
moduli space $H(X,L,r,e)$
parameterizing isomorphism classes of stable pairs $(E,\Phi)$, $E$
is a locally free sheaf on $X$ of rank $r$ and degree $e$ and $\Phi:E\to E\otimes_X L$
is a morphism of sheaves.

There is a torus action $\IC^\times \times H(X,L,r,e)\to
H(X,L,r,e)$, $t\times (E,\Phi) \to (E, t\Phi)$. The fixed points of the torus action
are stable pairs of the form
\be\label{eq:fixedpairsA}
E \simeq \oplus_{i=0}^nE_i, \qquad \Phi=\oplus_{i=0}^{n-1} {\varphi_i}
\ee
where $\varphi_i: E_i \to E_{i+1}\otimes_X L$, $i=0,\ldots, n-1$, all other components being trivial.
Note that the direct summand $E_i$, $0\leq i\leq n$, corresponds to the $\IC^\times$ character $t\to t^{-i}$. If $n=0$, $\Phi=0$, and $E=E_0$ must be stable bundle on $X$.

For $(r,e)\in \IZ_{\geq 1}\times \IZ$ coprime the torus fixed locus is smooth.
Given any connected component $\Xi$ of the fixed locus, the normal bundle
to $\Xi$ is isomorphic to the moving part of the tangent bundle to $H(X,L,r,e)$
restricted to $\Xi$,
\[
N_\Xi \simeq T^m_{H(X,L,r,e)}|_{\Xi}.
\]
Moreover, $N_\Xi$ decomposes in a direct sum of the form
\[
N_\Xi \simeq N_\Xi^+ \oplus N_\Xi^-
\]
where $N_\Xi^\pm$ is the direct sum of all $\IC^\times$ eigensheaves with
positive, respectively negative eigenvalues.
By definition, the index of the component $\Xi$ is $r^-_\Xi=r(N_\Xi^-)$.

The deformation theory of a Hitchin pair $(E,\Phi)$ is determined by the hypercohomology
of the two term complex on $X$
\be\label{eq:defthA}
0\to Hom_X(E,E) {\buildrel d\over \longto} Hom_X(E,E\otimes_X L) \to 0
\ee
where $d(f) = \Phi\circ f - f\otimes 1_L \circ \Phi$. If $(E,\Phi)$ is fixed by the torus action, the equivariant version of \eqref{eq:defthA} is
\be\label{eq:defthB}
0\to Hom_X(E,E) {\buildrel d\over \longto} Q\otimes Hom_X(E,E\otimes_X L) \to 0
\ee
where $E$ is of the form \eqref{eq:fixedpairsA} and $Q$ is the irreducible representation
of $\IC^\times$ with character $t\to t$. If $(E,\Phi)$ is a stable pair with $(r,e)$ coprime,
the 0-th hypercohomology group
of \eqref{eq:fixedpairsA} is isomorphic to $\IC$ while the 2nd hypercohomology group
vanishes. The 1st hypercohomology group is isomorphic to  the tangent space
$T_{[(E,\Phi)]}H(X,L,r,e)$.

The localization computations of the Hodge polynomial
of the moduli space of stable pairs of types $(2,1)$ and $(3,1)$ are reviewed
below.

\subsection{Rank $r=2$}
Let $e=1$. Then a fixed pair is either of the form $(E,0)$, with $E$ a stable
bundle of type $(2,1)$ on $X$ or
\[
E=E_0\oplus E_1,\qquad  \Phi=\left[\begin{array}{cc} 0 & 0\\ \varphi_0 & 0\\
\end{array}\right]
\]
with $E_0, E_1$ line bundles of degrees $e_0,e_1$, $e_0+e_1=1$ and
$\varphi_0:E_0\to E_1\otimes_X L$ a nonzero morphism.
In the second case the
 stability condition is equivalent to $e_1 \leq e_0$, while $\varphi\neq 0$
implies $e_0\leq e_1+2g-2+p$. Therefore
\[
{1\over 2}\leq e_0\leq g+{p-1\over 2}.
\]
In conclusion the fixed locus is a union of the form
\[
H(X,L,2,1)^{\IC^\times} \simeq M(2,1) \cup \bigcup_{e_0=1}^{g+[(p-1)/2]}
J_{e_0}(X) \times S^{2g-1+p-2e_0}(X)
\]
where $M(2,1)$ denotes the moduli space of stable bundles of type $(2,1)$ on $X$,  a smooth projective variety.
An elementary computation shows that $M(2,1)$ has index $0$ while each component
$J_{e_0}(X) \times S^{2g-1+p-2e_0}(X)$ has index $2e_0+g-2$, independent
of $p$. Then the Hodge polynomial of the moduli space of Hitchin pairs is
\be\label{eq:ranktwohodgeA}
\bal
& H_{(u,v)}(H(X,L,2,1))
= H_{(u,v)}(M(2,1))\\
&  + \sum_{e_0=1}^{g+[(p-1)/2]}
u^{2e_0+g-2} v^{2e_0+g-2} (1+u)^g(1+v)^gH_{(u,v)}(S^{2g-1+p-2e_0}(X))\\
\eal
\ee
Moreover, according to \cite{earl-kirwan,motive-bundles,hodge-twotwo}, the Hodge polynomial of $M(2,1)$ is
\be\label{eq:ranktwohodgeB}
H_{(u,v)}(M(2,1)) = (1+u)^g(1+v)^g{(1+u^2v)^g(1+uv^2)^g-(uv)^g(1+u)^g(1+v)^g\over
(1-uv)(1-u^2v^2)}
\ee
while the generating function for the Hodge polynomial of symmetric products is
\be\label{eq:ranktwohodgeC}
\sum_{n=0}^\infty x^n H_{(u,v)}(S^n(X)) = {(1+xu)^g(1+xv)^g \over (1-x)(1-xuv)}.
\ee
Repeating the computations of \cite{hitchin-selfd} in the present context yields
\be\label{eq:ranktwohdgeD}
\bal
& H_{(u,v)}(H(X,L,2,1)) =\\
& (1+u)^g(1+v)^g \bigg[
{(1+u^2v)^g(1+uv^2)^g\over
(1-uv)(1-u^2v^2)} +{(-1)^{p+1}\over 4}(uv)^{2g-2+p}{(1-u)^g(1-v)^g\over 1+uv} \\
& +(uv)^{2g-2+p} {(1+u)^g(1+v)^g\over 2(1-uv)}\bigg({g\over 1+u}+{g\over 1+v}
-{1\over 1-uv} -(2g-2+p)-{1\over 2}\bigg)\bigg]
\eal
\ee
\subsection{Rank $r=3$}
In this case the computation of the Poincar\'e polynomial has been done in \cite{Bettinumbers}.
Let $e=1$ for concretness; $e=2$ is analogous. The classification of fixed loci is more involved. There are four types of components.

$I)$ $(E,\Phi)=(E,0)$ with $E$ a rank 3 bundle on $X$ of degree $e=1$. This component is isomorphic to the moduli space $M(3,1)$
of stable rank 3 bundles on $X$ of degree $e=1$, which is a smooth projective variety.
Moreover, it has index 0 and according to \cite{earl-kirwan,hodge-rankthree}
\be\label{eq:rankthreeA}
\bal
& H_{(u,v)}(M(3,1)) = {(1+u)^g(1+v)^g\over (1-uv)(1-u^2v^2)^2(1-u^3v^3)}\\
& \bigg[
(1+u^2v^3)^g(1+u^3v^2)^g(1+uv^2)^g(1+u^2v)^g\\
& \ -(uv)^{2g-1}(1+uv)^2(1+u)^g(1+v)^g(1+uv^2)^g(1+u^2v)^g\\
&\  + (uv)^{3g-1}(1+uv+u^2v^2)(1+u)^{2g}(1+v)^{2g}\bigg]
\eal
\ee

$II)$  $E=E_0\oplus E_1$, $\Phi =\left[\begin{array}{cc} 0 & 0\\ \varphi_0 & 0\\
\end{array}\right]$, $E_0$ a degree $e_0$ line bundle,
$E_1$ a bundle of type  $(2,e_1)$ on $X$, and $\varphi_0:E_0\to E_1\otimes_X L$
a nontrivial morphism. Obviously $e_0+e_1=e$. The stability condition is equivalent to
the following two conditions
\begin{itemize}
\item ${e/3} \leq e_0 \leq {e/3}+ g-1+p/2$, which for $e=1$ yields
$1\leq e_0\leq g+p/2-2/3$.
\item The data
$(E_1\otimes L\otimes E_0^{-1}, {\varphi}_0\otimes 1_{E_0^{-1}})$
is a $\sigma$-stable Thaddeus pair of type $(2,e-3e_0+2(2g-2+p))=(2,1-3e_0+2(2g-2+p))$ (no fixed determinant) where
$\sigma =e_0/2 -e/6=e_0/2 -1/6$.
\end{itemize}
The index equals $3e_0-e+2g-2=3e_0+2g-3$, independent of $p$.
Therefore, repeating the computation in \cite[Sect. 4]{verlinde},
it follows that the contribution of fixed loci of this type to the Hodge
polynomial is
\be\label{eq:rankthreeB}
\bal
& H_{(u,v)}(\Xi^{II}(e_0)) = (uv)^{3e_0+2g-3} {(1+u)^{2g}(1+v)^{2g}\over 1-uv} \\
& \mathrm{Coeff}_{x^i} \bigg[\bigg({(uv)^{e_0+g}\over xu^2v^2-1}-
{(uv)^{2g-1-2e_0+p}\over x-uv}\bigg){(1+xu)^g(1+xv)^g\over (1-x)(1-xuv)}\bigg]
\eal
\ee
where $i=-2e_0+2g-2+p$.

$III)$ $E=E_0\oplus E_1$, $\Phi =\left[\begin{array}{cc} 0 & 0\\ \varphi_0 & 0\\
\end{array}\right]$, $E_0$ a bundle of type $(2,e_0)$ and $E_1$ a degree $e_1$
line bundle on $X$, $e_0+e_1=1$. In this case the dual pair $(E^\vee, \Phi^\vee\otimes
1_L)$ is a stable fixed pair of type $(II)$ with numerical invariants $(3,-e)=(3,-1)$.
Such fixed loci are labeled by an integer ${\bar e_0}$, $0\leq {\bar e}_0 \leq
g+p/2-4/3$ and their index is $3{\bar e}_0 +2g -1$.
Therefore their contribution to the Hodge polynomial is
\be\label{eq:rankthreeC}
\bal
& H_{(u,v)}(\Xi^{III}({\bar e}_0)) = (uv)^{3{\bar e}_0+2g-1}
{(1+u)^{2g}(1+v)^{2g}\over 1-uv} \\
& \mathrm{Coeff}_{x^i} \bigg[\bigg({(uv)^{{\bar e}_0+g}\over xu^2v^2-1}-
{(uv)^{2g-2-2e_0+p}\over x-uv}\bigg){(1+xu)^g(1+xv)^g\over (1-x)(1-xuv)}\bigg]
\eal
\ee
where $i=-2{\bar e}_0+2g-3+p$.

$IV)$  $E=E_0\oplus E_1\oplus E_2$, $\Phi =\left[\begin{array}{ccc} 0 & 0 & 0 \\ \varphi_0 & 0 & 0 \\ 0 & \varphi_1 & 0 \\
\end{array}\right]$, $E_0, E_1,E_2$ line bundles of degrees $e_0,e_1,e_2$ on $X$, $e_0+e_1+e_2=e=1$, and $\varphi_0:E_0\to E_1\otimes_X L$, $\varphi_1:E_1\to E_2\otimes_X L$
nontrivial morphisms. Let $m_1=e_1-e_0+2g-2+p$, $m_2=e_2-e_0+2g2+p$.
Then the stability conditions are equivalent to
\[
m_1,m_2\geq 0,\qquad m_1+2m_2\leq 3(2g-2+p), \qquad 2m_1+m_2\leq 3(2g-2+p).
\]
In addition, the following constraint holds by construction
\[
m_1+2m_2\equiv -e\ (\mathrm{mod}\ 3).
\]
Fixed loci of this type are isomorphic to a direct product of the form
$J_{e_0}(X)\times S^{m_1}(X)\times S^{m_2}(X)$.
The index is $8g-8+3p-m_1-m_2$.
Therefore their contribution to the Hodge polynomial is
\be\label{eq:rankthreeD}
\bal
& H_{(u,v)}(\Xi^{IV}_{(m_1,m_2)} )
= (uv)^{8g-8+3p-m_1-m_2}(1+u)^g(1+v)^g \\
& \mathrm{Coeff}_{x^{m_1}}\left({(1+xu)^g(1+xv)^g\over (1-x)(1-xuv)}\right)
\mathrm{Coeff}_{x^{m_2}}\left({(1+xu)^g(1+xv)^g\over (1-x)(1-xuv)}\right)\\
 \eal
\ee
In conclusion
\be
\bal
& H_{(u,v)}(H(X,L,3,1)) = H_{(u,v)}(M(3,1))
+ \sum_{e_0=1}^{g+[p/2-2/3]} H_{(u,v)}(\Xi^{II}(e_0)) \\ & +
\sum_{{\bar e}_0=0}^{g+[p/2-4/3]} H_{(u,v)}(\Xi^{III}({\bar e}_0)) +
\sum_{\substack{m_1,m_2\geq 0\\ 2m_1+m_2\leq 6g-6+3p\\ m_1+2m_2\leq
6g-6+3p\\ m_1+2m_2 \equiv 2\ (3)\\}} H_{(u,v)}(\Xi^{IV}_{(m_1,m_2)} )\eal
\ee
\subsection{Hausel-Rodriguez-Villegas Formula}
This subsection is a brief summary of the formulas of Hausel and Rodriguez-Villegas
\cite{HRV}, \cite{Mirror-Hodge}
for the Poincar\'e, respectively
Hodge polynomial of the moduli space $H(X,K_X,r,e)$ with
$(r,e)\in \IZ_{\geq}\times \IZ$ coprime.
Construct the following formal series
\[\CZ(q,x,y,T)=1+\sum_{k\geq1}T^kA_k(q,x,y)=1+\sum_{k\geq1}T^k\left(\sum_{|Y|=k}A_Y(q,x,y)\right)\]
where:
\[A_Y(q,x,y)=\prod_{z\in Y}\frac{(qxy)^{l(z)(2-2g)}(1+q^{h(z)}y^{l(z)}x^{l(z)+1})^g(1+q^{h(z)}x^{l(z)}y^{l(z)+1})^g}{(1-q^{h(z)}(xy)^{l(z)+1})(1-q^{h(z)}(xy)^{l(z)})}\]
where for $z=(i,j)\in Y$:
\[a(z)=Y_i-j,~l(z)=Y^t_j-i,~h(z)=a(z)+l(z)+1\] Define $H_r(q,x,y)$ in terms of the following recursive formula:

\[\sum_{r\geq 1}\sum_{k\geq1}H_r(q^k,-(-x)^k,-(-y)^k)B_r(q^k,-(-x)^k,-(-y)^k)\frac{T^{kr}}{k}=\log{\CZ(q,x,y,T)}\]
by comparing the coefficient of $T^{nk}$, where:
\[
B_r(q,x,y)=\frac{(qxy)^{(1-g)r(r-1)}(1+qx)^g(1+qy)^g}{(1-qxy)(1-q)}\]
Then
\be\label{eq:epoly}
 E_r(u,v)=H_r(1,u,v)
 \ee
 is conjectured in \cite{Mirror-Hodge} to be Hodge polynomial of the
moduli space $H(X,K_X,r,e)$.

\section{Recursion results for Hodge polynomials}\label{hodgeapp}
This section is basically a list of results for the Hodge polynomials of the moduli spaces
of Hitchin pairs determined by the double refinement the recursion formula
\eqref{eq:higgsrecursionA} and conjecture (\ref{refasympinvII}).
According to conjecture (\ref{refasympinvII}) the building blocks of asymptotic refined ADHM invariants are

\[
\bal
 \Omega^{(p)}_{\tableau{1}}
 = (-1)^p(uv)^{(1-g)/2} {(1-\lambda (uv^{-1})^{1/2})^g(1-\lambda(u^{-1}v)^{1/2})^g
 \over (1-\lambda (uv)^{1/2})(1-\lambda(uv)^{-1/2})}
 \eal
 \]
 \[
 \bal
 \Omega^{(p)}_{\tableau{2}}(\lambda,u,v)=&
 (uv)^{1-g-p/2}\lambda^{-p+2-2g}
 G(\lambda^2(uv)^{1/2},(uv)^{1/2},(uv^{-1})^{1/2})\\
 & G(\lambda,(uv)^{1/2},(uv^{-1})^{1/2})\\
  \eal
 \]
 \[
 \bal
 \Omega^{(p)}_{\tableau{1 1}}(\lambda,u,v)=& (uv)^{-p/2}\lambda^p
 G(\lambda^2(uv)^{-1/2},(uv)^{1/2},(uv^{-1})^{1/2})\\
 & G(\lambda, (uv)^{1/2},(uv^{-1})^{1/2})\\
   \eal
 \]

\[\bal \Omega^{(p)}_{\tableau{3}}(\lambda,u,v)=&(-1)^p(\lambda(uv)^{1/2})^{-3p+6(1-g)}G(\lambda^3uv,(uv)^{1/2},(uv^{-1})^{1/2})\\&G(\lambda^2(uv)^{1/2},(uv)^{1/2},(uv^{-1})^{1/2})G(\lambda,(uv)^{1/2},(uv^{-1})^{1/2})\eal\]

\[\bal \Omega^{(p)}_{\tableau{1 1 1}}(\lambda,u,v)=&(-1)^p(\lambda^{-1}(uv)^{1/2})^{-3p}G(\lambda^3(uv)^{-1},(uv)^{1/2},(uv^{-1})^{1/2})\\&G(\lambda^2(uv)^{-1/2},(uv)^{1/2},(uv^{-1})^{1/2})G(\lambda,(uv)^{1/2},(uv^{-1})^{1/2})\eal\]

\[ \Omega^{(p)}_{\tableau{2 1}}(\lambda,u,v)=(-1)^p(uv)^{-p}(\lambda (uv)^{1/2})^{2-2g}G(\lambda^3,(uv)^{1/2},(uv^{-1})^{1/2})G(\lambda,(uv)^{1/2},(uv^{-1})^{1/2})^2\]
where
\[
G(q,z,w)=
z^{(1-g)}\frac{(1-qw)^{g}(1-qw^{-1})^g}{(1-qz)(1-qz^{-1})}.
\]

  \bigskip
\framebox[1.7cm][l]{$1) \ {g=2}$} \par
\bigskip

\[\wH^{(0)}(2,1)(u,v)
=\frac{(1-u)^2(1-v)^2}{(\sqrt{uv})^9}(1+uv)(1 - 2 u^2 v - 2 u v^2 + 2 u^2 v^2)\]

\[\bal \wH^{(0)}(2,0)(u,v)
=& \frac{(1-u)^2(1-v)^2}{2(\sqrt{uv})^9(1+uv)}\\
&(2 + 4 u v - 4 u^2 v - 4 u v^2 + 7 u^2 v^2 - 6 u^3 v^2 + u^4 v^2 -
  6 u^2 v^3 + 12 u^3 v^3 \\&- 2 u^4 v^3 + u^2 v^4 - 2 u^3 v^4 +
  5 u^4 v^4)\eal \]

  \[\bal \wH^{(0)}(3,1)(u,v)=& \frac{(1-u)^2(1-v)^2}{(\sqrt{uv})^{19}}\\
  & (1 + u v - 2 u^2 v - 2 u v^2 + 3 u^2 v^2 - 4 u^3 v^2 + u^4 v^2 -
 4 u^2 v^3 + 8 u^3 v^3 \\&- 8 u^4 v^3 + 5 u^5 v^3 + u^2 v^4 -
 8 u^3 v^4 + 19 u^4 v^4 - 16 u^5 v^4 + 8 u^6 v^4 - 2 u^7 v^4 \\&+
 5 u^3 v^5 - 16 u^4 v^5 + 32 u^5 v^5 - 30 u^6 v^5 + 12 u^7 v^5 -
 2 u^8 v^5 + 8 u^4 v^6 - 30 u^5 v^6 \\&+ 43 u^6 v^6 - 32 u^7 v^6 +
 8 u^8 v^6 - 2 u^4 v^7 + 12 u^5 v^7 - 32 u^6 v^7 + 32 u^7 v^7 -
 12 u^8 v^7 \\&- 2 u^5 v^8 + 8 u^6 v^8 - 12 u^7 v^8 + 6 u^8 v^8)\eal\]

 \[ \bal \wH^{(1)}(2,1)(u,v)=& \frac{(1-u)^2(1-v)^2}{(\sqrt{uv})^{11}}\\
 &(1 + u v - 2 u^2 v - 2 u v^2 + 2 u^2 v^2 - 2 u^3 v^2 + u^4 v^2 -
 2 u^2 v^3 + 6 u^3 v^3 \\& - 2 u^4 v^3 + u^2 v^4 - 2 u^3 v^4 + 2 u^4
 v^4)\eal\]

\[\bal
\wH^{(1)}(2,0)(u,v)=& \frac{(1-u)^2(1-v)^2}{2(\sqrt{uv})^{11}(1+uv)}\\
&(2 + 4
u v - 4 u^2 v - 4 u v^2 + 6 u^2 v^2 - 8 u^3 v^2 + 2 u^4 v^2 -
 8 u^2 v^3 \\&+ 15 u^3 v^3 - 10 u^4 v^3 + u^5 v^3 + 2 u^2 v^4 -
 10 u^3 v^4 + 12 u^4 v^4 - 6 u^5 v^4 + u^3 v^5 \\&- 6 u^4 v^5 + 3 u^5
 v^5)\eal\]

 \[ \bal \wH^{(1)}(3,1)(u,v)=& \frac{(1-u)^2(1-v)^2}{(\sqrt{uv})^{25}}\\ &(6 v^{11} u^{11}-18 v^{10} u^{11}+18 v^9
   u^{11}-8 v^8 u^{11}+v^7 u^{11}-18 v^{11} u^{10}+60 v^{10} u^{10}\\&-76 v^9 u^{10}+42 v^8 u^{10}-10 v^7 u^{10}+v^6 u^{10}+18 v^{11} u^9-76 v^{10} u^9+121 v^9
   u^9\\&-98 v^8 u^9+42 v^7 u^9-8 v^6 u^9-8 v^{11} u^8+42 v^{10} u^8-98 v^9 u^8+124 v^8 u^8-84 v^7 u^8\\&+29 v^6 u^8-4 v^5 u^8+v^{11} u^7-10 v^{10} u^7+42 v^9 u^7-84
   v^8 u^7+91 v^7 u^7-56 v^6 u^7\\&+17 v^5 u^7-2 v^4 u^7+v^{10} u^6-8 v^9 u^6+29 v^8 u^6-56 v^7 u^6+58 v^6 u^6-34 v^5 u^6\\&+8 v^4 u^6-4 v^8 u^5+17 v^7 u^5-34 v^6
   u^5+33 v^5 u^5-16 v^4 u^5+5 v^3 u^5-2 v^7 u^4\\&+8 v^6 u^4-16 v^5 u^4+19 v^4 u^4-8 v^3 u^4+v^2 u^4+5 v^5 u^3-8 v^4 u^3+8 v^3 u^3\\&-4 v^2 u^3+v^4 u^2-4 v^3 u^2+3
   v^2 u^2-2 v u^2-2 v^2 u+v u+1)\eal\]

   \[\bal \wH^{(2)}(2,1)(u,v)=& \frac{(1-u)^2(1-v)^2}{(\sqrt{uv})^{13}}\\
   &(1 + u v) (1 - 2 u^2 v - 2 u v^2 + 2 u^2 v^2 + u^4 v^2 + 4 u^3 v^3 -
   4 u^4 v^3 \\&+ u^2 v^4 - 4 u^3 v^4 + 2 u^4 v^4)\eal\]

   \[\bal
   \wH^{(2)}(2,0)(u,v)=& \frac{(1-u)^2(1-v)^2}{2(\sqrt{uv})^{13}(1+uv)}\\
   &(2 + 4 u v - 4 u^2 v - 4 u v^2 + 6 u^2 v^2 - 8 u^3 v^2 + 2 u^4 v^2 -
 8 u^2 v^3 + 16 u^3 v^3 \\&- 12 u^4 v^3 + 4 u^5 v^3 + 2 u^2 v^4 -
 12 u^3 v^4 + 25 u^4 v^4 - 14 u^5 v^4 + 3 u^6 v^4 + 4 u^3 v^5 \\&-
 14 u^4 v^5 + 20 u^5 v^5 - 6 u^6 v^5 + 3 u^4 v^6 - 6 u^5 v^6 +
 5 u^6 v^6)\eal\]

 \[\bal\wH^{(2)}(3,1)(u,v)=& \frac{(1-u)^2(1-v)^2}{(\sqrt{uv})^{31}}\\
 &(10 v^{14} u^{14}-32 v^{13} u^{14}+36 v^{12}
   u^{14}-18 v^{11} u^{14}+3 v^{10} u^{14}-32 v^{14} u^{13}\\&+112 v^{13} u^{13}-154 v^{12} u^{13}+101 v^{11} u^{13}-32 v^{10} u^{13}+4 v^9 u^{13}+36 v^{14}
   u^{12}\\&-154 v^{13} u^{12}+269 v^{12} u^{12}-248 v^{11} u^{12}+123 v^{10} u^{12}-30 v^9 u^{12}+3 v^8 u^{12}\\&-18 v^{14} u^{11}+101 v^{13} u^{11}-248 v^{12}
   u^{11}+329 v^{11} u^{11}-250 v^{10} u^{11}+109 v^9 u^{11}\\&-24 v^8 u^{11}+v^7 u^{11}+3 v^{14} u^{10}-32 v^{13} u^{10}+123 v^{12} u^{10}-250 v^{11} u^{10}\\&+298
   v^{10} u^{10}-208 v^9 u^{10}+79 v^8 u^{10}-14 v^7 u^{10}+v^6 u^{10}+4 v^{13} u^9\\&-30 v^{12} u^9+109 v^{11} u^9-208 v^{10} u^9+226 v^9 u^9-146 v^8 u^9+53 v^7
   u^9\\&-8 v^6 u^9+3 v^{12} u^8-24 v^{11} u^8+79 v^{10} u^8-146 v^9 u^8+159 v^8 u^8\\&-96 v^7 u^8+30 v^6 u^8-4 v^5 u^8+v^{11} u^7-14 v^{10} u^7+53 v^9 u^7-96 v^8
   u^7\\&+98 v^7 u^7-58 v^6 u^7+17 v^5 u^7-2 v^4 u^7+v^{10} u^6-8 v^9 u^6+30 v^8 u^6\\&-58 v^7 u^6+59 v^6 u^6-34 v^5 u^6+8 v^4 u^6-4 v^8 u^5+17 v^7 u^5-34 v^6 u^5\\&+33
   v^5 u^5-16 v^4 u^5+5 v^3 u^5-2 v^7 u^4+8 v^6 u^4-16 v^5 u^4+19 v^4 u^4\\&-8 v^3 u^4+v^2 u^4+5 v^5 u^3-8 v^4 u^3+8 v^3 u^3-4 v^2 u^3+v^4 u^2-4 v^3 u^2\\&+3 v^2
   u^2-2 v u^2-2 v^2 u+v u+1)\eal\]

     \bigskip
\framebox[1.7cm][l]{$2) \ {g=3}$} \par
\bigskip

\[\bal\wH^{(0)}(2,1)(u,v)=& \frac{(1-u)^3(1-v)^3}{(\sqrt{uv})^{17}}\\
&(1 + u v - 3 u^2 v - 3 u v^2 + 2 u^2 v^2 - 3 u^3 v^2 + 3 u^4 v^2 -
 3 u^2 v^3 \\&+ 11 u^3 v^3 - 6 u^4 v^3 + 3 u^5 v^3 - u^6 v^3 +
 3 u^2 v^4 - 6 u^3 v^4 + 12 u^4 v^4 - 15 u^5 v^4 \\&+ 3 u^6 v^4 +
 3 u^3 v^5 - 15 u^4 v^5 + 12 u^5 v^5 - 6 u^6 v^5 - u^3 v^6 +
 3 u^4 v^6 \\&- 6 u^5 v^6 + 3 u^6 v^6)\eal\]

\[\bal\wH^{(0)}(2,0)(u,v)=& \frac{(1-u)^3(1-v)^3}{2(\sqrt{uv})^{17}(1+uv)}\\
&(2 + 4 u v - 6 u^2 v - 6 u v^2 + 6 u^2 v^2 - 12 u^3 v^2 + 6 u^4 v^2 -
 12 u^2 v^3 \\&+ 26 u^3 v^3 - 18 u^4 v^3 + 12 u^5 v^3 - 2 u^6 v^3 +
 6 u^2 v^4 - 18 u^3 v^4 + 47 u^4 v^4 \\&- 39 u^5 v^4 + 15 u^6 v^4 -
 u^7 v^4 + 12 u^3 v^5 - 39 u^4 v^5 + 57 u^5 v^5 - 33 u^6 v^5 \\&+
 9 u^7 v^5 - 2 u^3 v^6 + 15 u^4 v^6 - 33 u^5 v^6 + 39 u^6 v^6 -
 9 u^7 v^6 - u^4 v^7 \\&+ 9 u^5 v^7 - 9 u^6 v^7 + 7 u^7 v^7)\eal\]

\[\bal\wH^{(0)}(3,1)(u,v)=& \frac{(1-u)^3(1-v)^3}{(\sqrt{uv})^{37}}\\
&(15 v^{16} u^{16}-60 v^{15} u^{16}+102 v^{14}
   u^{16}-93 v^{13} u^{16}+45 v^{12} u^{16}-12 v^{11} u^{16}\\&+v^{10} u^{16}-60 v^{16} u^{15}+276 v^{15} u^{15}-537 v^{14} u^{15}+549 v^{13} u^{15}-324 v^{12}
   u^{15}\\&+102 v^{11} u^{15}-15 v^{10} u^{15}+v^9 u^{15}+102 v^{16} u^{14}-537 v^{15} u^{14}+1167 v^{14} u^{14}\\&-1407 v^{13} u^{14}+990 v^{12} u^{14}-417 v^{11}
   u^{14}+102 v^{10} u^{14}-12 v^9 u^{14}-93 v^{16} u^{13}\\&+549 v^{15} u^{13}-1407 v^{14} u^{13}+2003 v^{13} u^{13}-1776 v^{12} u^{13}+1020 v^{11} u^{13}\\&-362
   v^{10} u^{13}+72 v^9 u^{13}-6 v^8 u^{13}+45 v^{16} u^{12}-324 v^{15} u^{12}+990 v^{14} u^{12}\\&-1776 v^{13} u^{12}+2069 v^{12} u^{12}-1587 v^{11} u^{12}+798
   v^{10} u^{12}-251 v^9 u^{12}+42 v^8 u^{12}\\&-3 v^7 u^{12}-12 v^{16} u^{11}+102 v^{15} u^{11}-417 v^{14} u^{11}+1020 v^{13} u^{11}-1587 v^{12} u^{11}\\&+1659
   v^{11} u^{11}-1173 v^{10} u^{11}+537 v^9 u^{11}-156 v^8 u^{11}+21 v^7 u^{11}+v^{16} u^{10}\\&-15 v^{15} u^{10}+102 v^{14} u^{10}-362 v^{13} u^{10}+798 v^{12}
   u^{10}-1173 v^{11} u^{10}+1151 v^{10} u^{10}\\&-777 v^9 u^{10}+330 v^8 u^{10}-80 v^7 u^{10}+12 v^6 u^{10}+v^{15} u^9-12 v^{14} u^9\\&+72 v^{13} u^9-251 v^{12}
   u^9+537 v^{11} u^9-777 v^{10} u^9+731 v^9 u^9-456 v^8 u^9\\&+195 v^7 u^9-41 v^6 u^9+3 v^5 u^9-6 v^{13} u^8+42 v^{12} u^8-156 v^{11} u^8\\&+330 v^{10} u^8-456 v^9
   u^8+447 v^8 u^8-261 v^7 u^8+99 v^6 u^8-21 v^5 u^8\\&-3 v^{12} u^7+21 v^{11} u^7-80 v^{10} u^7+195 v^9 u^7-261 v^8 u^7+242 v^7 u^7\\&-147 v^6 u^7+45 v^5 u^7-10 v^4
   u^7+12 v^{10} u^6-41 v^9 u^6+99 v^8 u^6-147 v^7 u^6\\&+122 v^6 u^6-78 v^5 u^6+21 v^4 u^6-v^3 u^6+3 v^9 u^5-21 v^8 u^5+45 v^7 u^5\\&-78 v^6 u^5+63 v^5 u^5-30 v^4
   u^5+12 v^3 u^5-10 v^7 u^4+21 v^6 u^4-30 v^5 u^4\\&+34 v^4 u^4-12 v^3 u^4+3 v^2 u^4-v^6 u^3+12 v^5 u^3-12 v^4 u^3+13 v^3 u^3\\&-6 v^2 u^3+3 v^4 u^2-6 v^3 u^2+3 v^2
   u^2-3 v u^2-3 v^2 u+v u+1)\eal\]

\[\bal\wH^{(1)}(2,1)(u,v)=& \frac{(1-u)^3(1-v)^3}{(\sqrt{uv})^{19}}\\
&(1 + u v) (1 - 3 u^2 v - 3 u v^2 + 2 u^2 v^2 + 3 u^4 v^2 +
   9 u^3 v^3 - 6 u^4 v^3 \\&- u^6 v^3 + 3 u^2 v^4 - 6 u^3 v^4 +
   3 u^4 v^4 - 9 u^5 v^4 + 6 u^6 v^4 - 9 u^4 v^5 + 18 u^5 v^5 \\&-
   6 u^6 v^5 - u^3 v^6 + 6 u^4 v^6 - 6 u^5 v^6 + 3 u^6 v^6)\eal\]

   \[\bal\wH^{(1)}(2,0)(u,v)=& \frac{(1-u)^3(1-v)^3}{2(\sqrt{uv})^{19}(1+uv)}\\
   &(2 + 4 u v - 6 u^2 v - 6 u v^2 + 6 u^2 v^2 - 12 u^3 v^2 + 6 u^4 v^2 -
 12 u^2 v^3 \\&+ 26 u^3 v^3 - 18 u^4 v^3 + 12 u^5 v^3 - 2 u^6 v^3 +
 6 u^2 v^4 - 18 u^3 v^4 + 46 u^4 v^4 - 42 u^5 v^4 \\&+ 18 u^6 v^4 -
 4 u^7 v^4 + 12 u^3 v^5 - 42 u^4 v^5 + 65 u^5 v^5 - 63 u^6 v^5 +
 21 u^7 v^5 - 3 u^8 v^5 \\&- 2 u^3 v^6 + 18 u^4 v^6 - 63 u^5 v^6 +
 75 u^6 v^6 - 51 u^7 v^6 + 9 u^8 v^6 - 4 u^4 v^7 + 21 u^5 v^7 \\&-
 51 u^6 v^7 + 39 u^7 v^7 - 15 u^8 v^7 - 3 u^5 v^8 + 9 u^6 v^8 -
 15 u^7 v^8 + 5 u^8 v^8)\eal\]

\[\bal\wH^{(2)}(2,1)(u,v)=& \frac{(1-u)^3(1-v)^3}{(\sqrt{uv})^{21}}\\
&(1 + u v - 3 u^2 v - 3 u v^2 + 2 u^2 v^2 - 3 u^3 v^2 + 3 u^4 v^2 -
 3 u^2 v^3 + 11 u^3 v^3 \\&- 6 u^4 v^3 + 3 u^5 v^3 - u^6 v^3 +
 3 u^2 v^4 - 6 u^3 v^4 + 12 u^4 v^4 - 15 u^5 v^4 + 6 u^6 v^4 \\&-
 u^7 v^4 + 3 u^3 v^5 - 15 u^4 v^5 + 21 u^5 v^5 - 18 u^6 v^5 +
 9 u^7 v^5 - 2 u^8 v^5 - u^3 v^6 \\&+ 6 u^4 v^6 - 18 u^5 v^6 +
 30 u^6 v^6 - 24 u^7 v^6 + 6 u^8 v^6 - u^4 v^7 + 9 u^5 v^7 -
 24 u^6 v^7 \\&+ 21 u^7 v^7 - 9 u^8 v^7 - 2 u^5 v^8 + 6 u^6 v^8 -
 9 u^7 v^8 + 3 u^8 v^8)\eal\]

\[\bal\wH^{(2)}(2,0)(u,v)=& \frac{(1-u)^3(1-v)^3}{2(\sqrt{uv})^{21}(1+uv)}\\
&(2 + 4 u v - 6 u^2 v - 6 u v^2 + 6 u^2 v^2 - 12 u^3 v^2 + 6 u^4 v^2 -
 12 u^2 v^3 + 26 u^3 v^3 \\&- 18 u^4 v^3 + 12 u^5 v^3 - 2 u^6 v^3 +
 6 u^2 v^4 - 18 u^3 v^4 + 46 u^4 v^4 - 42 u^5 v^4 + 18 u^6 v^4 \\&-
 4 u^7 v^4 + 12 u^3 v^5 - 42 u^4 v^5 + 66 u^5 v^5 - 66 u^6 v^5 +
 30 u^7 v^5 - 6 u^8 v^5 - 2 u^3 v^6 \\&+ 18 u^4 v^6 - 66 u^5 v^6 +
 103 u^6 v^6 - 81 u^7 v^6 + 33 u^8 v^6 - 3 u^9 v^6 - 4 u^4 v^7 +
 30 u^5 v^7 \\&- 81 u^6 v^7 + 111 u^7 v^7 - 57 u^8 v^7 + 15 u^9 v^7 -
 6 u^5 v^8 + 33 u^6 v^8 - 57 u^7 v^8\\& + 57 u^8 v^8 - 15 u^9 v^8 -
 3 u^6 v^9 + 15 u^7 v^9 - 15 u^8 v^9 + 7 u^9 v^9)\eal\]

\bibliography{adhmref.bib}

\begin{thebibliography}{10}

\bibitem{motive-higgs}
L.~Alvarez-Consul, O.~Garcia-Prada, H.~J., M.~Logares, and A.~Schmitt.
\newblock The motive of the moduli space of rank 4 higgs bundles.
\newblock to appear.

\bibitem{YM-surface}
M.~F. Atiyah and R.~Bott.
\newblock The {Y}ang-{M}ills equations over {R}iemann surfaces.
\newblock {\em Philos. Trans. Roy. Soc. London Ser. A}, 308(1505):523--615,
  1983.

\bibitem{motivic-degzero}
K.~Behrend, J.~Bryan, and B.~Szendroi.
\newblock {Motivic degree zero Donaldson-Thomas invariants}.
\newblock arXiv:0909.5088.

\bibitem{infpairs}
I.~Biswas and S.~Ramanan.
\newblock An infinitesimal study of the moduli of {H}itchin pairs.
\newblock {\em J. London Math. Soc. (2)}, 49(2):219--231, 1994.

\bibitem{CV-I}
S.~Cecotti and C.~Vafa.
\newblock {{B}{P}{S} Wall Crossing and Topological Strings}.
\newblock hep-th/0910.2615.

\bibitem{higherrank}
W.-y. Chuang, D.~Diaconescu, and G.~Pan.
\newblock Rank two {A}{D}{H}{M} invariants and wallcrossing.
\newblock arXiv:10020579.

\bibitem{chamberII}
W.-y. Chuang, D.-E. Diaconescu, and G.~Pan.
\newblock Chamber structure and wallcrossing in the {A}{D}{H}{M} theory of
  curves {II}.
\newblock arXiv:0908.1119.

\bibitem{motive-bundles}
S.~Del~Ba{\~n}o.
\newblock On the motive of moduli spaces of rank two vector bundles over a
  curve.
\newblock {\em Compositio Math.}, 131(1):1--30, 2002.

\bibitem{DM-split}
F.~Denef and G.~W. Moore.
\newblock {Split states, entropy enigmas, holes and halos}.
\newblock arXi.org:hep-th/0702146.

\bibitem{Poincarepol_stabbundles}
U.~V. Desale and S.~Ramanan.
\newblock Poincar\'e polynomials of the variety of stable bundles.
\newblock {\em Math. Ann.}, 216(3):233--244, 1975.

\bibitem{chamberI}
D.-E. Diaconescu.
\newblock Chamber structure and wallcrossing in the {A}{D}{H}{M} theory of
  curves {I}.
\newblock arXiv:0904.4451.

\bibitem{modADHM}
D.~E. Diaconescu.
\newblock Moduli of {A}{D}{H}{M} sheaves and local {D}onaldson-{T}homas theory.
\newblock arXiv.org:0801.0820.

\bibitem{ruled}
D.-E. Diaconescu, B.~Florea, and N.~Saulina.
\newblock {A vertex formalism for local ruled surfaces}.
\newblock {\em Commun. Math. Phys.}, 265:201--226, 2006.

\bibitem{DM-crossing}
D.~E. Diaconescu and G.~W. Moore.
\newblock {Crossing the Wall: Branes vs. Bundles}.
\newblock arXiv.org:hep-th/0706.3193.

\bibitem{DG}
T.~Dimofte and S.~Gukov.
\newblock {Refined, Motivic, and Quantum}.
\newblock {\em Lett. Math. Phys.}, 91:1, 2010.

\bibitem{DGS}
T.~Dimofte, S.~Gukov, and Y.~Soibelman.
\newblock {Quantum Wall Crossing in N=2 Gauge Theories}.
\newblock arXiv:0912.1346.

\bibitem{earl-kirwan}
R.~Earl and F.~Kirwan.
\newblock The {H}odge numbers of the moduli spaces of vector bundles over a
  {R}iemann surface.
\newblock {\em Q. J. Math.}, 51(4):465--483, 2000.

\bibitem{EK-I}
T.~Eguchi and H.~Kanno.
\newblock {Five-dimensional gauge theories and local mirror symmetry}.
\newblock {\em Nucl. Phys.}, B586:331--345, 2000.

\bibitem{EK-II}
T.~Eguchi and H.~Kanno.
\newblock {Topological strings and Nekrasov's formulas}.
\newblock {\em JHEP}, 12:006, 2003.

\bibitem{laminations}
D.~Gaiotto, G.~Moore, and A.~Neitzke.
\newblock {Framed {B}{P}{S} states}.
\newblock arXiv:1006.0146.

\bibitem{betti-parabolic}
O.~Garc{\'{\i}}a-Prada, P.~B. Gothen, and V.~Mu{\~n}oz.
\newblock Betti numbers of the moduli space of rank 3 parabolic {H}iggs
  bundles.
\newblock {\em Mem. Amer. Math. Soc.}, 187(879):viii+80, 2007.

\bibitem{BPS-Higgs}
A.~Gholampour and A.~Sheshmani.
\newblock {{B}{P}{S} states via Higgs type bundles}.
\newblock in preparation.

\bibitem{Bettinumbers}
P.~B. Gothen.
\newblock The {B}etti numbers of the moduli space of stable rank {$3$} {H}iggs
  bundles on a {R}iemann surface.
\newblock {\em Internat. J. Math.}, 5(6):861--875, 1994.

\bibitem{hodgenumbers}
L.~G{\"o}ttsche.
\newblock Change of polarization and {H}odge numbers of moduli spaces of
  torsion free sheaves on surfaces.
\newblock {\em Math. Z.}, 223(2):247--260, 1996.

\bibitem{cohgroups}
G.~Harder and M.~S. Narasimhan.
\newblock On the cohomology groups of moduli spaces of vector bundles on
  curves.
\newblock {\em Math. Ann.}, 212:215--248, 1974/75.

\bibitem{Mirror-Hodge}
T.~Hausel.
\newblock Mirror symmetry and {L}anglands duality in the non-abelian {H}odge
  theory of a curve.
\newblock In {\em Geometric methods in algebra and number theory}, volume 235
  of {\em Progr. Math.}, pages 193--217. Birkh\"auser Boston, Boston, MA, 2005.

\bibitem{HRV}
T.~Hausel and F.~Rodriguez-Villegas.
\newblock Mixed {H}odge polynomials of character varieties.
\newblock {\em Invent. Math.}, 174(3):555--624, 2008.
\newblock With an appendix by Nicholas M. Katz.

\bibitem{hitchin-selfd}
N.~J. Hitchin.
\newblock The self-duality equations on a {R}iemann surface.
\newblock {\em Proc. London Math. Soc. (3)}, 55(1):59--126, 1987.

\bibitem{Hollowood:2003cv}
T.~J. Hollowood, A.~Iqbal, and C.~Vafa.
\newblock {Matrix Models, Geometric Engineering and Elliptic Genera}.
\newblock {\em JHEP}, 03:069, 2008.

\bibitem{IKP-I}
A.~Iqbal and A.-K. Kashani-Poor.
\newblock {Instanton counting and Chern-Simons theory}.
\newblock {\em Adv. Theor. Math. Phys.}, 7:457--497, 2004.

\bibitem{IKP-II}
A.~Iqbal and A.-K. Kashani-Poor.
\newblock {SU(N) geometries and topological string amplitudes}.
\newblock {\em Adv. Theor. Math. Phys.}, 10:1--32, 2006.

\bibitem{IKV}
A.~Iqbal, C.~Kozcaz, and C.~Vafa.
\newblock {The refined topological vertex}.
\newblock {\em JHEP}, 10:069, 2009.

\bibitem{J-I}
D.~Joyce.
\newblock Configurations in abelian categories. {I}. {B}asic properties and
  moduli stacks.
\newblock {\em Adv. Math.}, 203(1):194--255, 2006.

\bibitem{J-II}
D.~Joyce.
\newblock Configurations in abelian categories. {II}. {R}ingel-{H}all algebras.
\newblock {\em Adv. Math.}, 210(2):635--706, 2007.

\bibitem{J-III}
D.~Joyce.
\newblock Configurations in abelian categories. {III}. {S}tability conditions
  and identities.
\newblock {\em Adv. Math.}, 215(1):153--219, 2007.

\bibitem{J-IV}
D.~Joyce.
\newblock Configurations in abelian categories. {IV}. {I}nvariants and changing
  stability conditions.
\newblock {\em Adv. Math.}, 217(1):125--204, 2008.

\bibitem{genDTI}
D.~Joyce and Y.~Song.
\newblock A theory of generalized {D}onaldson-{T}homas invariants.
\newblock arxiv.org:0810.5645.

\bibitem{Katz:1996ht}
S.~H. Katz, D.~R. Morrison, and M.~Ronen~Plesser.
\newblock {Enhanced Gauge Symmetry in Type II String Theory}.
\newblock {\em Nucl. Phys.}, B477:105--140, 1996.

\bibitem{Konishi-I}
Y.~Konishi.
\newblock {Topological strings, instantons and asymptotic forms of
  Gopakumar-Vafa invariants}.
\newblock hep-th/0312090.

\bibitem{KS-review}
M.~Kontsevich and Y.~Soibelman.
\newblock Motivic {D}onaldson-{T}homas invariants: summary of results.
\newblock arXiv:0910:4315.

\bibitem{wallcrossing}
M.~Kontsevich and Y.~Soibelman.
\newblock Stability structures, {D}onaldson-{T}homas invariants and cluster
  transformations.
\newblock arXiv.org:0811.2435.

\bibitem{Lawrence:1997jr}
A.~E. Lawrence and N.~Nekrasov.
\newblock Instanton sums and five-dimensional gauge theories.
\newblock {\em Nucl. Phys.}, B513:239--265, 1998.

\bibitem{LLZ}
J.~Li, K.~Liu, and J.~Zhou.
\newblock Topological string partition functions as equivariant indices.
\newblock {\em Asian J. Math.}, 10(1):81--114, 2006.

\bibitem{adhm-rec}
S.~Mozgovoy.
\newblock {Solution of the motivic ADHM recursion formula}.
\newblock in preparation.

\bibitem{hodge-rankthree}
V.~Mu{\~n}oz.
\newblock Hodge polynomials of the moduli spaces of rank 3 pairs.
\newblock {\em Geom. Dedicata}, 136:17--46, 2008.

\bibitem{hodge-pairs}
V.~Mu{\~n}oz, D.~Ortega, and M.-J. V{\'a}zquez-Gallo.
\newblock Hodge polynomials of the moduli spaces of pairs.
\newblock {\em Internat. J. Math.}, 18(6):695--721, 2007.

\bibitem{hodge-twotwo}
V.~Mu{\~n}oz, D.~Ortega, and M.-J. V{\'a}zquez-Gallo.
\newblock Hodge polynomials of the moduli spaces of triples of rank {$(2,2)$}.
\newblock {\em Q. J. Math.}, 60(2):235--272, 2009.

\bibitem{N-refined}
K.~Nagao.
\newblock Refined open non-commutative {D}onaldson-{T}homas invariants for
  small crepant resolutions.
\newblock arxiv.org:0907.3784.

\bibitem{hilblect}
H.~Nakajima.
\newblock {\em Lectures on {H}ilbert schemes of points on surfaces}, volume~18
  of {\em University Lecture Series}.
\newblock American Mathematical Society, Providence, RI, 1999.

\bibitem{instcountA}
H.~Nakajima and K.~Yoshioka.
\newblock Instanton counting on blowup. {I}. 4-dimensional pure gauge theory.
\newblock {\em Invent. Math.}, 162(2):313--355, 2005.

\bibitem{instcountB}
H.~Nakajima and K.~Yoshioka.
\newblock Instanton counting on blowup. {II}. {$K$}-theoretic partition
  function.
\newblock {\em Transform. Groups}, 10(3-4):489--519, 2005.

\bibitem{Nekrasov:2002qd}
N.~A. Nekrasov.
\newblock {Seiberg-Witten Prepotential From Instanton Counting}.
\newblock {\em Adv. Theor. Math. Phys.}, 7:831--864, 2004.

\bibitem{MR1085642}
N.~Nitsure.
\newblock Moduli space of semistable pairs on a curve.
\newblock {\em Proc. London Math. Soc. (3)}, 62(2):275--300, 1991.

\bibitem{projective}
A.~Schmitt.
\newblock Projective moduli for {H}itchin pairs.
\newblock {\em Internat. J. Math.}, 9(1):107--118, 1998.

\bibitem{simpsonI}
C.~T. Simpson.
\newblock Moduli of representations of the fundamental group of a smooth
  projective variety. {I}.
\newblock {\em Inst. Hautes \'Etudes Sci. Publ. Math.}, (79):47--129, 1994.

\bibitem{simpsonII}
C.~T. Simpson.
\newblock Moduli of representations of the fundamental group of a smooth
  projective variety. {II}.
\newblock {\em Inst. Hautes \'Etudes Sci. Publ. Math.}, (80):5--79 (1995),
  1994.

\bibitem{Tachikawa:2004ur}
Y.~Tachikawa.
\newblock {Five-dimensional Chern-Simons terms and Nekrasov's instanton
  counting}.
\newblock {\em JHEP}, 02:050, 2004.

\bibitem{verlinde}
M.~Thaddeus.
\newblock Stable pairs, linear systems and the {V}erlinde formula.
\newblock {\em Invent. Math.}, 117(2):317--353, 1994.

\bibitem{ranktwo}
Y.~Toda.
\newblock On a computation of rank two {D}onaldson-{T}homas invariants.
\newblock arxiv.org:0912.2507.

\bibitem{KY-Betti}
K.~Yoshioka.
\newblock The {B}etti numbers of the moduli space of stable sheaves of rank
  {$2$} on a ruled surface.
\newblock {\em Math. Ann.}, 302(3):519--540, 1995.

\bibitem{chamber}
K.~Yoshioka.
\newblock Chamber structure of polarizations and the moduli of stable sheaves
  on a ruled surface.
\newblock {\em Internat. J. Math.}, 7(3):411--431, 1996.

\end{thebibliography}
 \bibliographystyle{abbrv}
\end{document}